\documentclass[12pt]{article}
\usepackage{amsthm, amsmath, amssymb, amsfonts, url, booktabs, tikz, setspace, fancyhdr, bm}
\usepackage{hyperref}
\usepackage{geometry}
\usepackage{hyperref, enumerate}
\usepackage[shortlabels]{enumitem}
\usepackage[babel]{microtype}
\usepackage[english]{babel}
\usepackage[capitalise]{cleveref}
\usepackage{multirow}
\usepackage{pifont}

\usepackage{algorithm}
\usepackage{algorithmic}
\usepackage{bbm,tkz-graph,subcaption}
\usepackage{csquotes}
\usepackage{mathrsfs}
\usepackage{mathabx}
\usetikzlibrary{patterns}
\usetikzlibrary{shapes}
\usetikzlibrary{calc}
\usepackage{bbm}
\usepackage{makecell}

\usepackage[utf8]{inputenc} 
\usepackage[T1]{fontenc}
\usepackage{amsfonts}
\usepackage{amscd}
\usepackage{graphicx}
\usepackage{enumitem}
\usepackage{verbatim}
\usepackage{hyperref}
\usepackage{amsmath,caption}
\usepackage{url,pdfpages,xcolor,framed,color}
\usepackage{todonotes}

\newtheorem{thr}{Theorem}[section]

\newtheorem{p}[thr]{Problem}

\newtheorem{prop}[thr]{Proposition}
\newtheorem{conj}[thr]{Conjecture}
\theoremstyle{definition}

\newtheorem*{defi*}{Definition}

\newcommand{\NP}{\textsf{NP}}
\newcommand{\cmark}{\ding{51}}

\newcommand*{\myproofname}{Proof}

\title{Computer-assisted graph theory: a survey}

\date{}
\author{
{\sc Jorik JOOKEN\footnote{Department of Computer Science, KU Leuven Campus Kulak-Kortrijk, 8500 Kortrijk, Belgium}}
\footnote{E-mail address: jorik.jooken@kuleuven.be}}

\begin{document}
\maketitle
\begin{abstract}
Computers and algorithms play an ever-increasing role in obtaining new results in graph theory. In this survey, we present a broad range of techniques used in computer-assisted graph theory, including the exhaustive generation of all pairwise non-isomorphic graphs within a given class, the use of searchable databases containing graphs and invariants as well as other established and emerging algorithmic paradigms. We cover approaches based on mixed integer linear programming, semidefinite programming, dynamic programming, SAT solving, metaheuristics and machine learning. The techniques are illustrated with numerous detailed results covering several important subareas of graph theory such as extremal graph theory, graph coloring, structural graph theory, spectral graph theory, regular graphs, topological graph theory, special sets in graphs, algebraic graph theory and chemical graph theory. We also present some smaller new results that demonstrate how readily a computer-assisted graph theory approach can be applied once the appropriate tools have been developed.
\end{abstract}


\section{Introduction}
Computer-assisted graph theory is concerned with the development, implementation and execution of algorithms in order to gain insights into a problem in graph theory. These insights can take many forms in accordance with the many different use cases of computer-assisted graph theory. Typically, the main focus in the literature is on the graph-theoretical results rather than the algorithms even though the algorithms are often a very important ingredient to arrive at these results. As a consequence, the description of the graph-theoretical results is often much more detailed than the description of the algorithms. In some cases even only the graph-theoretical results are described without mentioning the use of an algorithm at all, since the output of the algorithm was part of the inspiration gathering process and this is often not explicitly written down in a paper. An example of this is finding a counterexample to a conjecture by means of a computer-assisted approach and subsequently providing a purely mathematical proof (without the need for a computer) that this is indeed a valid counterexample. The clear disadvantage of this practice is that it prevents other researchers from learning and following a similar computer-assisted approach for related problems and therefore slows down the potential progress that can be made. In this survey, we attempt to address this problem by giving an overview of several existing techniques in the broad domain of computer-assisted graph theory as well as their use cases. Throughout this survey, we will also give several references to publicly available software libraries that are very useful for researchers wishing to engage in this domain themselves and we will give plenty of detailed examples to illustrate the use cases.

The use of computers in graph theory has a rich and evolving history. The birth of the domain of computer-assisted graph theory coincides with Appel and Haken's proof of the Four Color Theorem~\cite{AH77,AHK77}, stating that the regions of any planar map (equivalently, the faces of any planar graph) can be colored with at most four colors so that no two adjacent regions share the same color. This was the first major result obtained with the help of a computer-assisted graph theory approach by reducing the problem to one where thousands of configurations have to be checked. However, this result was not without controversy for several reasons. At the time, the computer-assisted nature of the proof made people wonder about the question of what constitutes a proof. Moreover, several researchers pointed out some sources of errors (some of them being more severe than others). As a response to this, Appel and Haken later addressed and corrected the errors in a book~\cite{AH89} containing a 400-page supplement that aids the reader in checking the necessary cases. Since then, Robertson, Sanders, Seymour and Thomas~\cite{RSST97} presented a new computer-assisted proof that is shorter and also leads to a quadratic time algorithm for producing the sought coloring (instead of the quartic time algorithm derived from Appel and Haken's proof). This convinced nearly the entire community of the validity of the proof. In 2007, any remaining doubts were dispelled by a formally verifiable computer-assisted proof using the proof assistant \texttt{Rocq}~\cite{G07} (previously known as \texttt{Coq}). However, to this day there is still no classical mathematical human-made proof of the Four Color Theorem. Over the years, computer-assisted graph theory has developed into a mature research domain that led to several new exciting results that will be discussed later.

In this survey, we exclusively deal with simple undirected graphs, i.e., ordered pairs $(V,E)$ where $V$ is a set representing the vertices and $E \subseteq \{\{u,v\}~|~u,v \in V\text{ and }u \neq v\}$ is a set representing the edges (unordered pairs of distinct vertices). Henceforth, we will simply refer to these as \textit{graphs}. We stress that many ideas appearing in computer-assisted graph theory are also applicable to other combinatorial objects (with which one can often naturally associate a graph) such as digraphs, multigraphs, hypergraphs, matroids, combinatorial designs and objects stemming from algebra and finite geometry (see~\cite{CRSSM00,LTS89,MR08,MMM07,MW05,PRS96} for examples of computer-assisted methods playing an important role for other combinatorial objects). Throughout this survey we use standard graph theory terminology. We refer the reader to~\cite{D17} for any terminology not defined in the current survey.

The rest of this survey is structured as follows: In Section~\ref{sec:graphGeneration} we discuss algorithms that can be used to generate all graphs satisfying given constraints. We focus on datasets of graphs (graph censuses) that were generated using such algorithms and discuss techniques that are used to avoid generating isomorphic copies. In Section~\ref{sec:databases} we present several publicly available databases containing graphs and various invariants. Such databases can be used to search for specific graphs that one is interested in. In Section~\ref{sec:otherParadigms} we discuss other algorithmic paradigms that are used in a computer-assisted graph theory context, namely mixed integer linear programming, semidefinite programming, dynamic programming, SAT-based approaches, metaheuristics and approaches based on machine learning. Finally, in Section~\ref{sec:newResults}, we close this survey by presenting some new results that can be obtained by elementary means using a computer-assisted graph theory approach. 

\section{Graph generation}
\label{sec:graphGeneration}
For many problems in graph theory, it is useful to generate certain graph classes. Often, one is not interested in all of these graphs, but rather in a subset of these graphs that satisfy a given property. Therefore, these graph generators can often be reused if the graph class is a widely studied one, whereas the algorithm for checking the property might have to be implemented from scratch (depending on whether software for checking this property is already publicly available or not). Most algorithms for generating graphs expect at least one parameter as input, namely the order $n$ of the graphs to be generated, and output a list of all pairwise non-isomorphic graphs of order $n$ in the graph class that the algorithm intends to generate. 

There are plenty of reasons why one could be interested in generating these graphs. For example, it is useful to find counterexamples to conjectures of the form \lq Every graph satisfying property $P_1$ also satisfies property $P_2$\rq. There are several different approaches that one could take to find counterexamples to such a conjecture, each one having its own strengths and weaknesses. For instance, one could generate all graphs and check which ones satisfy $P_1$, but do not satisfy $P_2$, or one could generate all graphs satisfying $P_1$ and check which ones do not satisfy $P_2$, or maybe one could even directly generate all graphs satisfying $P_1$, but not satisfying $P_2$. Although these subtleties seem to be details, they have a big impact on the graph generation algorithms that should be developed and their speed. 

Apart from generating counterexamples, graph generation is also often useful to find suitable \textit{gadgets}, i.e., graphs that will be used as part of a larger construction. Clearly, if the gadget is not too large, finding a gadget might be much easier than finding the larger construction for a generation algorithm (as we will see later, for many classes of graphs the number of pairwise non-isomorphic graphs of order $n$ quickly increases for increasing $n$).

Graph generation is also useful for problems where one is interested in finding the smallest order $n$ of a graph that satisfies a certain property $P_1$. If the algorithm is an exhaustive algorithm (and hence generates all graphs of a given order satisfying $P_1$), checking every graph with order at most $n'$ and not finding any suitable graph leads to the lower bound $n'<n$. On the other hand, if the algorithm finds a suitable graph satisfying $P_1$, this clearly provides an upper bound for $n$. The latter situation where one is looking for upper bounds for $n$ often occurs when generating all graphs that satisfy $P_1$ is difficult, but generating all graphs that satisfy $P_1$ and also another property $P_2$ is easier. For example, it is often useful to encode some symmetry properties of the graphs using $P_2$, see for example~\cite{CJ25b}.

In some cases, graph generation algorithms can also be used to prove that there are only finitely many graphs satisfying a certain property. This is achieved by implementing a recursive graph generation algorithm that employs certain pruning rules. The algorithm could recurse indefinitely without terminating, but in case these pruning rules are powerful enough the algorithm could terminate and this serves as a proof that there are only finitely many graphs satisfying a certain property.

Graph generation often also complements classical proofs by being able to rule out certain small graphs. These small graphs are often cumbersome to handle manually and this leads to large case distinctions. Handling them with a computer-assisted approach on the other hand leads to proofs that are easier to follow. A nice example of this can be found in~\cite{K20}.

In this section, we first give an overview of several publicly available algorithms for generating certain graph classes and give references to several downloadable lists of graphs (graph censuses). Graphs are ubiquitous and appear in a lot of contexts such as social networks, collaboration networks, road maps, biomedical datasets, web graphs and neural networks (see for example~\cite{RA15}). However, in the current paper we limit ourselves to graph classes that have received considerable attention from a graph-theoretical point of view (based on for example their structure, properties and algorithmic tractability results) rather than from the point of view of network science~\cite{B13}. Moreover, we remark that most of the graph generators that we discuss are exhaustive. In other words, they are capable of generating all graphs within a certain class (sometimes satisfying additional properties that the user can specify). Therefore, algorithms that focus on random generation of graphs are also out of scope (we refer the interested reader to~\cite{Green21} for more information about this interesting topic). Finally, we close this section by discussing some of the underlying principles of these graph generators and by discussing further considerations.

\subsection{Generation algorithms and graph censuses}
\label{sec:generationAlgorithms}
The \texttt{nauty} package~\cite{MP14} is a software package containing several useful programs for researchers in computer-assisted graph theory (this package is regularly updated and we refer the interested reader to~\cite{NautyUserGuide} for a detailed description of the latest version). It also contains several well-known graph generators that we will now discuss. 

The \texttt{geng} generator allows one to generate all pairwise non-isomorphic graphs on $n$ vertices. Therefore, this graph generator is the most general one and applicable to almost any problem that involves the generation of graphs. However, as shown in Table~\ref{tab:smallGraphs} the number of pairwise non-isomorphic graphs on $n$ vertices grows very quickly as the order $n$ increases (see also~\cite{OEISSmallGraphs}) and therefore in many cases it is useful to have additional constraints (and specialized graph generators) that allow one to effectively generate larger graphs satisfying these constraints.

\begin{table}[h]
    \centering
    \begin{tabular}{|c| r |}
\hline
$n$ & \thead{Number of pairwise\\ non-isomorphic graphs\\ on $n$ vertices}\\
\hline
1 & 1\\
2 & 2\\
3 & 4\\
4 & 11\\
5 & 34\\
6 & 156\\
7 & 1044\\
8 & 12 346\\
9 & 274 668\\
10 & 12 005 168\\
11 & 1 018 997 864\\
12 & 165 091 172 592\\
13 & 50 502 031 367 952\\
\hline
\end{tabular}
    \caption{The number of pairwise non-isomorphic graphs on $n$ vertices.}
    \label{tab:smallGraphs}
\end{table}

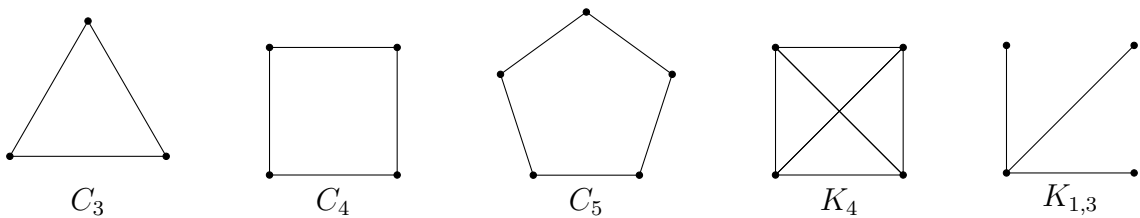
\begin{figure}[h!]
\centering
\begin{tikzpicture}[scale=0.4]
  \def\sides{3}
  \def\radius{3}

  \foreach \i in {1,2,3} {
    \draw ({360/\sides * (\i + 1)-30}:\radius) -- ({360/\sides * \i-30}:\radius);
  }
  
  \foreach \i in {1,2,3} {
    \draw[fill,black] ({360/\sides * \i-30}:\radius) circle (3pt);
  }
  \node at (0,-3) {$C_3$};
\end{tikzpicture}
\hspace{10mm}
\begin{tikzpicture}[scale=0.4]
  \def\sides{4}
  \def\radius{3}

  \foreach \i in {1,2,3,4} {
    \draw ({360/\sides * (\i + 1)+45}:\radius) -- ({360/\sides * \i+45}:\radius);
  }
  
  \foreach \i in {1,2,3,4} {
    \draw[fill,black] ({360/\sides * \i+45}:\radius) circle (3pt);
  }
  \node at (0,-3) {$C_4$};
\end{tikzpicture}
\hspace{10mm}
\begin{tikzpicture}[scale=0.4]
  \def\sides{5}
  \def\radius{3}

  \foreach \i in {1,2,3,4,5} {
    \draw ({360/\sides * (\i + 1)+18}:\radius) -- ({360/\sides * \i+18}:\radius);
  }
  
  \foreach \i in {1,2,3,4,5} {
    \draw[fill,black] ({360/\sides * \i+18}:\radius) circle (3pt);
  }
  \node at (0,-3.3) {$C_5$};
\end{tikzpicture}
\hspace{10mm}
\begin{tikzpicture}[scale=0.4]
  \def\sides{4}
  \def\radius{3}

  \foreach \i in {1,2,3,4} {
    \draw ({360/\sides * (\i + 1)+45}:\radius) -- ({360/\sides * \i+45}:\radius);
    \draw ({360/\sides * (\i + 2)+45}:\radius) -- ({360/\sides * \i+45}:\radius);
  }
  
  \foreach \i in {1,2,3,4} {
    \draw[fill,black] ({360/\sides * \i+45}:\radius) circle (3pt);
  }
  \node at (0,-3) {$K_4$};
\end{tikzpicture}
\hspace{10mm}
\begin{tikzpicture}[scale=0.4]
  \def\sides{4}
  \def\radius{3}


  \draw ({360/\sides * (1)+45}:\radius) -- ({360/\sides * 2+45}:\radius);
  \draw ({360/\sides * (0)+45}:\radius) -- ({360/\sides * 2+45}:\radius);
  \draw ({360/\sides * (3)+45}:\radius) -- ({360/\sides * 2+45}:\radius);
  
  \foreach \i in {1,2,3,4} {
    \draw[fill,black] ({360/\sides * \i+45}:\radius) circle (3pt);
  }
  \node at (0,-3) {$K_{1,3}$};
\end{tikzpicture}
\caption{The graphs $C_3$, $C_4$, $C_5$, $K_4$ and $K_{1,3}$ (the last one is also known as the \textit{claw graph).}}\label{fig:smallH}
\end{figure}

The \texttt{geng} algorithm allows one to specify further constraints on the number of edges and the minimum and maximum degree of the vertices in the graph. It also allows one to only generate (bi-)connected graphs, chordal graphs~\cite{B93}, split graphs~\cite{M03}, perfect graphs~\cite{CRST06,L72}, bipartite graphs and $H$-free graphs for every $H \in \{C_3,C_4,C_5,K_4,K_{1,3}\}$ (see Fig.~\ref{fig:smallH}). Here, we say that a graph $G$ is $H$-free if $G$ does not have an induced subgraph that is isomorphic to $H$.

Other graph generators in the \texttt{nauty} package include \texttt{gentreeg} (for generating all trees on $n$ vertices, where bounds for the diameter and an upper bound for the maximum degree can also be specified), \texttt{genktreeg} (for generating all $k$-trees on $n$ vertices, i.e., the maximal graphs with treewidth $k$~\cite{ND08}), \texttt{genquarticg} (for generating 4-regular graphs) and \texttt{genspecialg} (for generating a large number of special graphs that are easy to generate due to their structural description, such as generalized Petersen graphs or flower snarks~\cite{I75}).

Another graph class that attracted a lot of attention in the literature is the class of \textit{planar graphs}~\cite{NC88} (i.e., graphs that can be embedded in the plane without crossing edges). An example of an interesting open conjecture for planar graphs is the Matheson-Tarjan conjecture~\cite{MT96} stating that every \textit{planar triangulation} $G$ (i.e., maximal planar graph) of sufficiently large order $n$ has a \textit{dominating set} (i.e., a vertex set $D \subseteq V(G)$ such that every vertex in $V(G)\setminus D$ has a neighbor in $D$) of size at most $n/4$. The most well-known graph generator for generating planar graphs is \texttt{plantri}~\cite{BGGMTW05,BM05,BM07,VM24}. Apart from planar graphs, this generator also allows one to generate many important subclasses thereof such as planar triangulations and \textit{planar quadrangulations} (where every face has size 3 or 4, respectively), planar bipartite graphs, triangulations of a disk, planar cubic and quartic graphs and \textit{Apollonian triangulations} (i.e., the planar 3-trees). Moreover, \texttt{plantri} also allows one to specify further constraints on the vertex-connectivity, edge-connectivity, minimum degree and girth (the precise constraints depend on the subclass). Another important subclass worth mentioning is the class of \textit{fullerenes}~\cite{KHOCS85}, i.e., planar cubic graphs in which all faces have size 5 or 6. These graphs are important in chemistry as they model the atom-bond structure of chemical fullerenes (carbon allotropes) for the discovery of which Kroto, Curl, and Smalley received the Nobel Prize in Chemistry in 1996. These graphs can be efficiently generated using the graph generators \texttt{fullgen}~\cite{BD97} and \texttt{buckygen}~\cite{BGM12,GM15}. Table~\ref{tab:planarGraphs} shows the number of pairwise non-isomorphic graphs on $n$ vertices for several subclasses of planar graphs (see also~\cite{OEISPlanarTri,OEISPlanar,OEISPlanarQuad}) and supports our earlier claim that specialized generators are often desirable.

\begin{table}[h]
    \centering
    \resizebox{\textwidth}{!}{
    \begin{tabular}{|c| r | r | r |}
\hline
$n$ & \thead{Number of pairwise\\ non-isomorphic 3-connected \\
planar graphs on $n$ vertices} & \thead{Number of pairwise\\ non-isomorphic 3-connected \\
planar triangulations on $n$ vertices}  & \thead{Number of pairwise\\ non-isomorphic 3-connected \\
planar quadrangulations on $n$ vertices}\\
\hline
4 & 1 & 1 & 0\\
5 & 2 & 1 & 0\\
6 & 7 & 2 & 0\\
7 & 34 & 5 & 0\\
8 & 257 & 14 & 1\\
9 & 2606 & 50 & 0\\
10 & 32 300 & 233 & 1\\
11 & 440 564 & 1249 & 1\\
12 & 6 384 634 & 7595 & 3\\
13 & 96 262 938 & 49 566 & 3\\
14 & 1 496 225 352 & 339 722 & 11\\
15 & 23 833 988 129 & 2 406 841 & 18\\
16 & 387 591 510 244 & 17 490 241 & 58\\
17 & 6 415 851 530 241 & 129 664 753 & 139\\
18 & 107 854 282 197 058 & 977 526 957 & 451\\
\hline
\end{tabular}
}
    \caption{The number of pairwise non-isomorphic graphs on $n$ vertices for several subclasses of planar graphs.}
    \label{tab:planarGraphs}
\end{table}

A \textit{Hamiltonian cycle} in a graph is a cycle that contains every vertex of the graph precisely once. A \textit{$k$-regular graph} is a graph in which every vertex has $k$ neighbors. A longstanding conjecture by Sheehan states that every 4-regular graph containing a Hamiltonian cycle contains at least one more Hamiltonian cycle distinct from the first one~\cite{Sh75}. This is one of the many examples of conjectures and reasons making it interesting to generate regular graphs. These graphs can be effectively generated using \texttt{genreg}~\cite{Me99}. Moreover, \texttt{genreg} also allows one to specify a lower bound on the girth of the graphs to be generated.

Among the regular graphs, the 3-regular graphs (further referred to as \textit{cubic graphs}) received the most attention in the literature. Several of the most studied conjectures in graph theory such as the Petersen Coloring Conjecture~\cite{J88}, the 5-Flow Conjecture~\cite{T54}, the Berge-Fulkerson Conjecture~\cite{F71}, the Fan-Raspaud Conjecture~\cite{FR94} and the Cycle Double Cover Conjecture~\cite{S79,S73} are either directly stated for cubic graphs or it is known that a smallest counterexample (should one exist) must be cubic. Apart from being cubic, for many conjectures it is known that a counterexample cannot admit a proper 3-edge-coloring. We will refer to cubic graphs that do not admit a proper 3-edge-coloring as \textit{snarks}. The Petersen graph is without doubt the most famous example of a snark. We remark that in many cases, even much more is known about a potential smallest counterexample such as lower bounds on the girth or \textit{cyclic edge-connectivity}~\cite{H00,MM20}. We say that a graph $G$ is \textit{cyclically $k$-edge-connected} if $G$ does not contain a set of edges $C$ such that $G-C$ contains at least two components having at least one cycle and $|C|<k$. The cyclic edge-connectivity of a graph $G$ is then equal to the largest integer $k$ such that $G$ is cyclically $k$-edge-connected. Cubic graphs can be effectively generated using \texttt{minibaum}~\cite{B96} and \texttt{snarkhunter}~\cite{BGHM13,BGM11}. The former generator also allows one to only output bipartite graphs and specify a lower bound on the girth, whereas the latter generator allows one to specify a lower bound on the girth and only output snarks.

\begin{table}[h!]
    \centering
    \resizebox{\textwidth}{!}{
    \begin{tabular}{|c||c|c|c|c|c|c|c|}
\hline
Graph class & \thead{Degree\\sequence} & \thead{Girth} & \thead{Automorphisms} & \thead{2-Factor} & \thead{Distance} & \thead{Edge\\coloring} & \thead{Planar}\\
\hline
\multirow{2}{*}{\makecell{\textit{Small regular graphs}\\\textit{of large girth (cages)}~\cite{BMS95,E04,EJ12,EMMN11,MMN98,Me99}}} & \multirow{2}{*}{\cmark} & \multirow{2}{*}{\cmark} & \multirow{2}{*}{} & \multirow{2}{*}{} & \multirow{2}{*}{} & \multirow{2}{*}{} & \multirow{2}{*}{} \\  
& & & & & & & \\
\hline
\makecell{\textit{Small bi-regular graphs}\\\textit{of large girth}~\cite{GJV24}} & \cmark & \cmark & & & & & \\  
\hline
\makecell{\textit{Edge-girth-regular graphs}~\cite{GJ25}} & \cmark & \cmark & & & & & \\  
\hline
\makecell{\textit{Vertex-girth-regular graphs}~\cite{JJP24}} & \cmark & \cmark & & & & & \\  
\hline
\makecell{\textit{Regular graphs with girth $g$}\\\textit{but no $(g{+}1)$-cycles}~\cite{EJJ26}} & \cmark & \cmark & & & & & \\  
\hline
\multirow{2}{*}{\makecell{\textit{Strongly regular}\\\textit{graphs}~\cite{StronglyRegularGraphs,BV22,CDS06,HS01,StronglyRegularGraphs2,MS01,CombinatorialDataSpence,S00}}} & \multirow{2}{*}{\cmark} & \multirow{2}{*}{\cmark} & \multirow{2}{*}{} & \multirow{2}{*}{} & \multirow{2}{*}{\cmark} & \multirow{2}{*}{} & \multirow{2}{*}{} \\  
& & & & & & & \\
\hline
\makecell{\textit{Nanocones}~\cite{BV11}} & \cmark & \cmark & & & & & \cmark \\  
\hline
\makecell{Minimal Cayley graphs~\cite{GK25}} & \cmark & & \cmark & & & & \\  
\hline
\makecell{Vertex-transitive graphs~\cite{HR20}} & \cmark & & \cmark & & & & \\
\hline
\makecell{Edge-transitive graphs~\cite{EncyclopediaOfGraphs}} & \cmark & & \cmark & & & & \\
\hline
\multirow{2}{*}{\makecell{Cubic vertex-transitive and\\ Cayley graphs~\cite{PSV13,PSV15}}} & \multirow{2}{*}{\cmark} & \multirow{2}{*}{} & \multirow{2}{*}{\cmark} & \multirow{2}{*}{} & \multirow{2}{*}{} & \multirow{2}{*}{} & \multirow{2}{*}{} \\  
& & & & & & & \\
\hline
\makecell{Cubic edge-transitive graphs~\cite{CP25}} & \cmark & & \cmark & & & & \\  
\hline
\makecell{Cubic arc-transitive graphs~\cite{CombinatorialDataConder,CD02}} & \cmark & & \cmark & & & & \\  
\hline
\makecell{Cubic semisymmetric graphs~\cite{CombinatorialDataConder,CMMP06}} & \cmark & & \cmark & & & & \\  
\hline
\makecell{Quartic edge-transitive graphs~\cite{PW20}} & \cmark & & \cmark & & & & \\  
\hline
\multirow{2}{*}{\makecell{Quartic arc-transitive graphs~\cite{PSV13,PSV15}}} & \multirow{2}{*}{\cmark} & \multirow{2}{*}{} & \multirow{2}{*}{\cmark} & \multirow{2}{*}{} & \multirow{2}{*}{} & \multirow{2}{*}{} & \multirow{2}{*}{} \\  
& & & & & & & \\
\hline
\makecell{Quartic 2-arc-transitive graphs~\cite{P09}} & \cmark & & \cmark & & & & \\  
\hline
\makecell{Circulant graphs~\cite{CombinatorialDataMcKay}} & \cmark & & \cmark & & & & \\  
\hline
\multirow{2}{*}{\makecell{Many classes of graphs\\defined by automorphism group~\cite{PWebsite}}} & \multirow{2}{*}{\cmark} & \multirow{2}{*}{} & \multirow{2}{*}{\cmark} & \multirow{2}{*}{} & \multirow{2}{*}{} & \multirow{2}{*}{} & \multirow{2}{*}{} \\  
& & & & & & & \\
\hline
\multirow{2}{*}{\makecell{\textit{Cycle permutation}\\\textit{graphs}~\cite{GR25}}} & \multirow{2}{*}{\cmark} & \multirow{2}{*}{} & \multirow{2}{*}{} & \multirow{2}{*}{\cmark} & \multirow{2}{*}{} & \multirow{2}{*}{} & \multirow{2}{*}{} \\  
& & & & & & & \\
\hline
\makecell{\textit{Permutation snarks}~\cite{GR25}} & \cmark & & & \cmark & & \cmark & \\  
\hline
\multirow{2}{*}{\makecell{Large graphs of given\\maximum degree and diameter~\cite{DH94,E01,LS08,MS12}}} & \multirow{2}{*}{\cmark} & \multirow{2}{*}{} & \multirow{2}{*}{} & \multirow{2}{*}{} & \multirow{2}{*}{\cmark} & \multirow{2}{*}{} & \multirow{2}{*}{} \\  
& & & & & & & \\
\hline
\makecell{\textit{Alternating plane graphs}~\cite{AHSSV15}} & \cmark & & & & & & \cmark \\  
\hline
\makecell{\textit{Several classes}\\ \textit{of planar graphs}~\cite{BFLPV10,MYS22}} & \cmark & & & & & & \cmark \\  
\hline
\makecell{Quartic integral graphs~\cite{MW15}} & \cmark & & & & & & \\  
\hline
\makecell{Highly irregular graphs~\cite{CombinatorialDataMcKay}} & \cmark & & & & & & \\  
\hline
\makecell{Halin graphs~\cite{TaylorWebsite}} &  & & & & & & \cmark \\  
\hline
\end{tabular}
}
    \caption{An overview of different graph generators or censuses and keywords (`Degree sequence', `Girth', `Automorphisms', `2-Factor', `Distance', `Edge coloring', `Planar') that are directly relevant for the definition of the graph class. Graph classes for which we found source code to generate graphs are indicated in italics.}
    \label{tab:classesOverview1}
\end{table}

\begin{table}[h!]
    \centering
    \resizebox{\textwidth}{!}{
    \begin{tabular}{|c||c|c|c|c|c|c|c|}
\hline
Graph class & \thead{Long\\cycles} & \thead{Criticality} & \thead{Spectrum} & \thead{Intersection} & \thead{Genus} & \thead{Vertex\\coloring} & \thead{Distance}\\
\hline
\makecell{\textit{(Almost) hypohamiltonian graphs}~\cite{AMW97,GZ17}} & \cmark & & & & & & \\
\hline
\makecell{\textit{$K_2$-(hypo)hamiltonian graphs}~\cite{GRWZ24,GRZ24}} & \cmark & & & & & & \\
\hline
\makecell{\textit{Perihamiltonian graphs}~\cite{FMTVZ21}} & \cmark & & & & & & \\
\hline
\makecell{\textit{Platypus graphs}~\cite{GNZ20}} & \cmark & & & & & & \\
\hline
\makecell{\textit{Graphs with few Hamiltonian cycles}~\cite{GMZ20}} & \cmark & & & & & & \\
\hline
\makecell{\textit{HIST-critical graphs}~\cite{GNRZ24}} & & \cmark & & & & & \\
\hline
\makecell{Minimal Ramsey graphs~\cite{AM24,GR13,R12}} & & \cmark & & & & & \\
\hline
\makecell{Torus obstructions~\cite{MW18}} & & \cmark & & & \cmark & & \\
\hline
\multirow{2}{*}{\makecell{\textit{Minimal $\mathcal{F}$-free obstructions}\\\textit{to $H$-coloring}~\cite{GJORS24,GS18,HMRSV15,XJGH25a,XJGH25b}}} & \multirow{2}{*}{} & \multirow{2}{*}{\cmark} & \multirow{2}{*}{} & \multirow{2}{*}{} & \multirow{2}{*}{} & \multirow{2}{*}{\cmark} & \multirow{2}{*}{} \\
& & & & & & & \\
\hline
\multirow{2}{*}{\makecell{\textit{(Maximal) triangle-free}\\\textit{($k$-chromatic) graphs}~\cite{BBH98,BBH00,BGS12,G20}}} & \multirow{2}{*}{} & \multirow{2}{*}{\cmark} & \multirow{2}{*}{} & \multirow{2}{*}{} & \multirow{2}{*}{} & \multirow{2}{*}{\cmark} & \multirow{2}{*}{} \\
& & & & & & & \\
\hline
\multirow{2}{*}{\makecell{Extremal graphs\\for Tur{\'a}n-type problems~\cite{A16}}} & \multirow{2}{*}{} & \multirow{2}{*}{\cmark} & \multirow{2}{*}{} & \multirow{2}{*}{} & \multirow{2}{*}{} & \multirow{2}{*}{} & \multirow{2}{*}{} \\
& & & & & & & \\
\hline
\multirow{2}{*}{\makecell{Extremal graphs\\for Zarankiewicz numbers~\cite{GHO00}}} & \multirow{2}{*}{} & \multirow{2}{*}{\cmark} & \multirow{2}{*}{} & \multirow{2}{*}{} & \multirow{2}{*}{} & \multirow{2}{*}{} & \multirow{2}{*}{} \\
& & & & & & & \\
\hline
\makecell{\textit{Nut graphs}~\cite{CFG18}} & & & \cmark & & & & \\
\hline
\makecell{Conduction-isomorphic graphs~\cite{BFGJ25}} & & & \cmark & & & & \\
\hline
\makecell{Cospectral graphs~\cite{HS04}} & & & \cmark & & & & \\
\hline
\multirow{3}{*}{\makecell{Graphs with maximal\\algebraic connectivity\\for a given order, diameter or girth~\cite{EKJS24}}} & & & & & & & \\
& & & \cmark & & & & \cmark \\
& & & & & & & \\
\hline
\makecell{Interval graphs~\cite{YSKU20}} & & & & \cmark & & & \\
\hline
\makecell{Permutation graphs~\cite{YSKU20}} & & & & \cmark & & & \\
\hline
\makecell{\textit{Cographs}~\cite{JPD18}} & & & & \cmark & & & \\
\hline
\makecell{\textit{Block graphs}~\cite{ESY24}} & & & & \cmark & & & \\
\hline
\makecell{\textit{Graphs embedded on a surface}~\cite{B25,S17}} & & & & & \cmark & & \\
\hline
\makecell{Distance-hereditary graphs~\cite{YQU21}} & & & & & & & \cmark \\
\hline
\makecell{Ptolemaic graphs~\cite{YQU21}} & & & & & & & \cmark \\
\hline

\end{tabular}
}
    \caption{An overview of different graph generators or censuses and keywords (`Long cycles', `Criticality', `Spectrum', `Intersection', `Genus', `Vertex coloring', `Distance') that are directly relevant for the definition of the graph class. Graph classes for which we found source code to generate graphs are indicated in italics.}
    \label{tab:classesOverview2}
\end{table}

We believe that the above graph generators are among the most well-known ones, but of course, there exist many other graph generators that are either more recent than the previously discussed ones, focus on a narrower class of graphs or for which the source code itself is not available for download, but instead the graph census obtained by running the graph generator is downloadable. Among others, we mention those shown in Table~\ref{tab:classesOverview1} and Table~\ref{tab:classesOverview2}. For each of the graph classes, we added keywords that are directly relevant for the definition of the graph class. For more details, we refer the interested reader to the corresponding papers for a definition of these graph classes and to~\cite{StronglyRegularGraphs,HoGMetaDirectory,RegularGraphsPage,PWebsite,CombinatorialDataWanless} for online overviews of such graph generators and censuses. Some graph classes received so much attention from the computer-assisted graph theory viewpoint that listing all relevant sources is impossible and therefore we point to relevant surveys containing such sources instead. Among others, we mention small regular graphs of large girth (cages)~\cite{EJ12}, large graphs of given maximum degree and diameter~\cite{MS12}, Ramsey graphs~\cite{R12} and strongly regular graphs~\cite{BV22}.

\subsection{Underlying principles}
The generated graphs are often built in an incremental fashion (e.g., by adding edges one by one), where backtracking is used to explore all possibilities (complemented by pruning rules to speed up the generation). An important source of inefficiency stems from the fact that building graphs in two different ways can lead to isomorphic graphs and the search space is duplicated if the search continues from these isomorphic copies. Therefore, it is important that the algorithms try to avoid generating isomorphic copies (both for intermediate graphs and the final graphs) as this leads to redundancy. In this subsection, we will discuss several techniques that can be used to achieve this.

Before we introduce these techniques, however, it is necessary to discuss three closely related problems.
\begin{itemize}
    \item The first problem is the \textit{graph isomorphism problem}, which asks one to decide for two given graphs $G$ and $H$ whether they are isomorphic. This problem is a classical example of a problem for which it is not known whether it is \NP-complete or solvable in polynomial time~\cite{KST12}. However, for many special cases (where there are restrictions on the input graphs $G$ and $H$), the problem can be solved in polynomial time~\cite{BGM82,B90,L82,M80,M04}. Moreover, in practice it is often easy to decide that $G$ and $H$ are not isomorphic using heuristics such as (variants of) the Weisfeiler-Lehman algorithm~\cite{WL68}.
    \item The second problem is the \textit{graph automorphism problem}, which asks one to compute (the generators of) the automorphism group of a given graph $G$. This problem is in fact polynomial-time equivalent to the graph isomorphism problem~\cite{M79}, but in practice it is often a bit harder.
    \item The third problem asks to compute a \textit{canonical labeling} of a given graph $G$. Here, two graphs $G$ and $H$ are isomorphic if and only if they have the same canonical labeling. This problem is the hardest of the three problems.
\end{itemize}
It is important to note that these three problems received a lot of attention in the literature and they can often be solved very efficiently for practical inputs arising from the graph generation algorithms. As a result, these graph generation algorithms often make calls to efficient existing software libraries that solve (some of) these three problems. Among others, we mention \texttt{nauty}~\cite{M81,MP14}, \texttt{Traces}~\cite{MP14,P08}, \texttt{bliss}~\cite{JK07,JK11},
\texttt{conauto}~\cite{LCF13,LF09}, \texttt{nishe}~\cite{TN08},
\texttt{saucy}~\cite{DLSM04,DSM08,KSM10} and \texttt{dejavu}~\cite{AS21a,AS21b,AS21c}. Moreover, in this context it is also appropriate to mention \texttt{GAP}~\cite{GAP} and \texttt{Magma}~\cite{Magma}, which are software libraries that are useful for many algebraic purposes and \texttt{SageMath}~\cite{sagemath}, which has built-in support for several graph generators. We note that several of these software libraries also allow us to compute the~\textit{vertex orbits} or the~\textit{edge orbits} of a graph. Here, two vertices (edges) of a graph $G$ are in the same vertex (edge) orbit if and only if there exists an automorphism of $G$ that maps the first vertex (edge) to the second vertex (edge).

\subsubsection{Techniques to avoid generating isomorphic graphs}

We now discuss five techniques that can be used to avoid generating isomorphic copies of a graph (the interested reader is referred to~\cite{B00,KS99} for a more extensive overview of this large research domain).

\begin{itemize}
    \item One of the simplest techniques is to keep track of a list of generated graphs and every time a new graph is generated, one can compare this graph with every graph that is already in the list. If the new graph is not isomorphic to any of the graphs that are already in the list, this graph is added to the list and otherwise it is discarded. In order to make this idea feasible, one has to be able to efficiently solve the graph isomorphism problem, since this problem has to be solved a linear number of times (in the length of the list) every time a new graph is added. The strength of this idea is that it is very simple and that the graph isomorphism problem is the easiest of the three mentioned problems to solve. The clear disadvantage, however, is that all graphs need to be stored in memory (even after they are written to the output) and this approach does not scale well when there are a lot of graphs involved due to memory limitations. Therefore, this idea is only recommended for situations where there are relatively few graphs and the sizes of the graphs are too large for other techniques to handle.

    We are in fact not aware of any well-known publicly available graph generators that use this first technique, because this technique is rarely the best option.
    
    \item For the second technique, one keeps track of both the generated graphs themselves and also their canonical forms. If these canonical forms are stored in an appropriate data structure such as a balanced binary search tree or a hash table, it is no longer necessary to compare a new element with every element in the data structure (containing $n$ elements), but only with $O(\log(n))$ elements (in the case of a balanced binary search tree) or $O(1)$ elements (while also computing a hash function in the case of a hash table)~\cite{K97}. For this technique, one needs to solve the third problem repeatedly by computing canonical forms. As with the first technique, the memory disadvantage remains. We remark, however, that for both of these techniques, if a small number of isomorphic graphs can be tolerated, it is also possible to only store a fixed number of graphs or canonical forms in memory by simply stopping to add new graphs or use a replacement strategy when the memory capacity is reached. In other words, one can always use less memory at the expense of longer computation times (and this is sometimes desirable, especially in for example a high-performance computing environment as we will briefly discuss later). Due to the simplicity of these two techniques, very little overhead is necessary and these techniques can be very fast in practice if the right conditions are met.

    The algorithms described in~\cite{EJJ26,GJ25} are examples of algorithms that use this technique.
    
    \item The third technique is the \textit{orderly generation} of graphs~\cite{F78,R78}, which was introduced independently by Farad\v{z}ev and Reed. This technique allows one to generate every graph exactly once (both intermediate and final graphs) without generating isomorphic copies and without storing all graphs (or canonical forms) in memory. Therefore, this technique works very well in situations where a lot of graphs have to be generated. The idea of this technique is to define a canonical labeling for every graph in every isomorphism class and only generate the graph with that canonical labeling for every isomorphism class (and no other graphs). Hence, the algorithm will never generate two isomorphic graphs. The canonical labeling should be chosen so that larger canonically labeled graphs can be constructed from smaller canonically labeled graphs. Clearly, the main challenge of designing such an algorithm consists of defining an appropriate canonical labeling such that, given an arbitrary labeling of one of the intermediate graphs, it is not too difficult to check that this labeling is the canonical one. This is highly dependent on the graph class that the algorithm should generate and typically this approach works best for graph classes where it is not too difficult to describe all automorphisms of the graphs (or the subgraphs of those graphs that the algorithm generates). Indeed, if one has an explicit description of all automorphisms of the graph, it suffices to generate all of these automorphisms and check whether the current labeling is the canonical one. 
    
    We remark that in some cases it is desirable to allow some isomorphic copies to be generated (by relaxing the canonical labeling constraint and therefore obtaining a computationally cheaper check to verify if a labeling should be allowed or not) and filter the isomorphic copies in a post-processing step. This is typically the case when it is expected that the algorithm will only output a few graphs.

    The \texttt{genreg} algorithm~\cite{Me99} and the algorithms described in~\cite{B96,GR25} are examples of algorithms that use this technique.

    \item The fourth technique that we will discuss is the \textit{canonical construction path method} introduced by McKay~\cite{M98}. Similarly as before, it is not necessary to store all graphs in memory using this technique. It is meant to be used for algorithms that recursively generate graphs by repeatedly expanding a smaller graph (for example by adding edges one by one) in all possible ways by means of backtracking. We will refer to such an operation as an \textit{expansion} (e.g., adding one edge) and refer to the reverse operation as a \textit{reduction} (e.g., removing one edge). In the remainder of this subsection, we will assume that the expansion operation consists of adding one edge to a graph, since this makes the explanation easy to follow, but we remark that the principle works for other expansions as well. There are three distinct situations that may lead to isomorphic copies (see Fig.~\ref{fig:situationOneAndTwo} for the first and second situation and Fig.~\ref{fig:situationThree} for the third situation).

    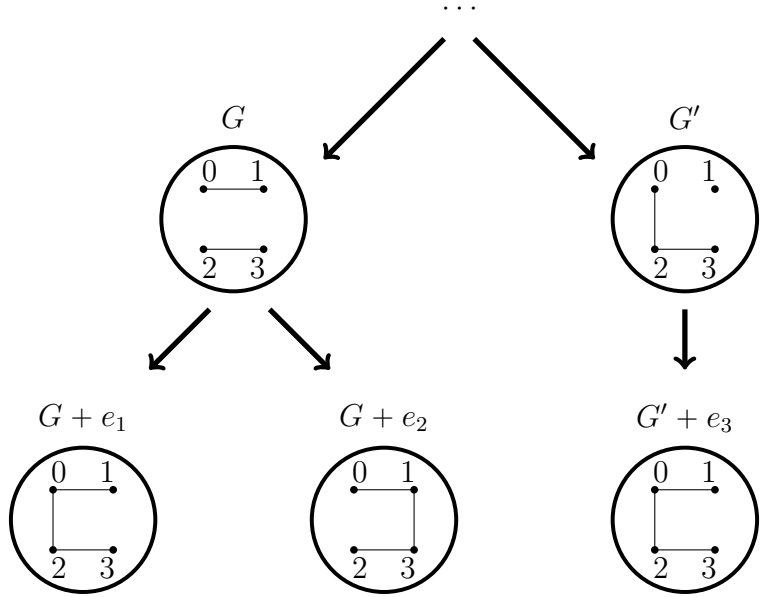
\begin{figure}[h!]
\begin{center}
\begin{tikzpicture}[scale=0.8]

\draw[->,line width=2pt] (0.1,4) -- (-0.9,3);
\draw[->,line width=2pt] (1.1,4) -- (2.1,3);
\draw[->,line width=2pt] (8,4) -- (8,3);

\draw[->,line width=2pt] (4,8.5) -- (2,6.5);
\draw[->,line width=2pt] (4.5,8.5) -- (6.5,6.5);
\node at (4.3,9) {$\cdots$};

\begin{scope}[shift={(7.5,5)}]
\draw (0,0) --(1,0);
\draw (0,0) --(0,1);
\fill (0,0) circle (1.8pt);
\fill (0,1) circle (1.8pt);
\fill (1,0) circle (1.8pt);
\fill (1,1) circle (1.8pt);

\node at (0.1, -0.3) {$2$};
\node at (0.9, -0.3) {$3$};
\node at (0.9, 1.3) {$1$};
\node at (0.1, 1.3) {$0$};

\draw[line width=1.5pt] (0.5,0.5) circle [radius=1.2];
\node at (0.5, 2.2) {$G'$};
\end{scope}

\begin{scope}[shift={(0,5)}]
\draw (0,0) --(1,0);
\draw (0,1) --(1,1);
\fill (0,0) circle (1.8pt);
\fill (0,1) circle (1.8pt);
\fill (1,0) circle (1.8pt);
\fill (1,1) circle (1.8pt);

\node at (0.1, -0.3) {$2$};
\node at (0.9, -0.3) {$3$};
\node at (0.9, 1.3) {$1$};
\node at (0.1, 1.3) {$0$};

\draw[line width=1.5pt] (0.5,0.5) circle [radius=1.2];
\node at (0.5, 2.2) {$G$};
\end{scope}

\begin{scope}[shift={(7.5,0)}]
\draw (0,0) --(1,0);
\draw (0,1) --(1,1);
\draw (0,0) --(0,1);
\fill (0,0) circle (1.8pt);
\fill (0,1) circle (1.8pt);
\fill (1,0) circle (1.8pt);
\fill (1,1) circle (1.8pt);

\node at (0.1, -0.3) {$2$};
\node at (0.9, -0.3) {$3$};
\node at (0.9, 1.3) {$1$};
\node at (0.1, 1.3) {$0$};

\draw[line width=1.5pt] (0.5,0.5) circle [radius=1.2];
\node at (0.5, 2.2) {$G'+e_3$};
\end{scope}

\begin{scope}[shift={(-2.5,0)}]
\draw (0,0) --(1,0);
\draw (0,1) --(1,1);
\draw (0,0) --(0,1);
\fill (0,0) circle (1.8pt);
\fill (0,1) circle (1.8pt);
\fill (1,0) circle (1.8pt);
\fill (1,1) circle (1.8pt);

\node at (0.1, -0.3) {$2$};
\node at (0.9, -0.3) {$3$};
\node at (0.9, 1.3) {$1$};
\node at (0.1, 1.3) {$0$};

\draw[line width=1.5pt] (0.5,0.5) circle [radius=1.2];
\node at (0.5, 2.2) {$G+e_1$};
\end{scope}

\begin{scope}[shift={(2.5,0)}]
\draw (0,0) --(1,0);
\draw (0,1) --(1,1);
\draw (1,0) --(1,1);
\fill (0,0) circle (1.8pt);
\fill (0,1) circle (1.8pt);
\fill (1,0) circle (1.8pt);
\fill (1,1) circle (1.8pt);

\node at (0.1, -0.3) {$2$};
\node at (0.9, -0.3) {$3$};
\node at (0.9, 1.3) {$1$};
\node at (0.1, 1.3) {$0$};

\draw[line width=1.5pt] (0.5,0.5) circle [radius=1.2];
\node at (0.5, 2.2) {$G+e_2$};
\end{scope}

\end{tikzpicture}
\end{center}
\caption{A part of the recursion tree showing the first and the second situation that may lead to isomorphic copies.}\label{fig:situationOneAndTwo} 
\end{figure}

     In the first situation, there exists a graph $G$ and distinct edges $e_1$ and $e_2$ in the complement of $G$ such that they are in the same edge orbit. Therefore, clearly $G+e_1$ will be isomorphic to $G+e_2$. This source of isomorphism can be prevented by calculating the edge orbits of the complement of $G$, and only consider one edge for each edge orbit. 
    
    In the second situation, however, there exists another graph $G'$ (not isomorphic to $G$) and edges $e_1$ (in the complement of $G$) and $e_3$ (in the complement of $G'$) such that the graph $G+e_1$ was already generated by the backtracking algorithm before and the graphs $G+e_1$ and $G'+e_3$ are isomorphic. 

    In the third situation, there exists a graph $G''$ and edges $e_4$ and $e_5$ in the complement of $G''$ such that $e_4$ and $e_5$ are not in the same edge orbit of the complement of $G''$, yet the graphs $G''+e_4$ and $G''+e_5$ are isomorphic (this is called \textit{pseudosimilarity}). An example of this is shown in Fig.~\ref{fig:situationThree}.
    \begin{figure}[ht]
        \centering
        \begin{tikzpicture}[scale=0.65]

    \draw[->,line width=2pt] (8.5,6.2) -- (10.5,4.5);
    \draw[->,line width=2pt] (3.4,6.2) -- (1.4,4.5);
    
  \begin{scope}[shift={(0,0)}]
    \foreach \x in {0,1,...,5}{
      \draw[fill] (\x*360/6:2.5) circle (1.8pt);
    }
    \draw (1*360/6:2.5)--({(1+1)*360/6}:2.5);
    \draw (2*360/6:2.5)--({(2+1)*360/6}:2.5);
    \draw (3*360/6:2.5)--({(3+1)*360/6}:2.5);
    \draw (5*360/6:2.5)--({(5+1)*360/6}:2.5);

    \draw (1*360/6:2.5)--({(0)*360/6+30}:3.5);
    \draw (0*360/6:2.5)--({(0)*360/6+30}:3.5);

    \draw (2*360/6:2.5)--({(2)*360/6+30}:3.5);
    \draw (3*360/6:2.5)--({(2)*360/6+30}:3.5);

    \draw (4*360/6:2.5)--({(4)*360/6+30}:3.5);
    \draw (5*360/6:2.5)--({(4)*360/6+30}:3.5);

    \draw (0*360/6:2.5)--(1*360/6:2.5);
    
    \foreach \x in {0,2,4}{
      \draw[fill] (\x*360/6+30:3.5) circle (1.8pt);
      \draw (\x*360/6:3.5)--({\x*360/6}:2.5);
    }
    \foreach \x in {0,2,4}{
      \draw[fill] (\x*360/6:3.5) circle (1.8pt);

      \node at (1.20, 2.6) {$0$};
      \node at (2.95, 2.3) {$1$};
      \node at (3.45, 0.4) {$2$};
      \node at (2.85, 0.4) {$3$};
      \node at (1.20, -1.75) {$4$};
      \node at (0, -2.75) {$5$};
      \node at (-2.0, -2.75) {$6$};
       
      \node at (-1.10, -1.75) {$7$};
      \node at (-2.95, 0.4) {$8$};
      \node at (-3.3, 2.3) {$9$};

      \node at (-1.30, 3.3) {$10$};
      \node at (-1.10, 2.6) {$11$};
      
      \node at (0, 4.8) {$G''+e_4$};
      
      \draw[line width=1.5pt] (-0.2,0) circle [radius=4.3];
    }
  \end{scope}

  \begin{scope}[shift={(12,0)}]
    \foreach \x in {0,1,...,5}{
      \draw[fill] (\x*360/6:2.5) circle (1.8pt);
    }
    \draw (1*360/6:2.5)--({(1+1)*360/6}:2.5);
    \draw (2*360/6:2.5)--({(2+1)*360/6}:2.5);
    \draw (3*360/6:2.5)--({(3+1)*360/6}:2.5);
    \draw (5*360/6:2.5)--({(5+1)*360/6}:2.5);

    \draw (1*360/6:2.5)--({(0)*360/6+30}:3.5);
    \draw (0*360/6:2.5)--({(0)*360/6+30}:3.5);

    \draw (2*360/6:2.5)--({(2)*360/6+30}:3.5);
    \draw (3*360/6:2.5)--({(2)*360/6+30}:3.5);

    \draw (4*360/6:2.5)--({(4)*360/6+30}:3.5);
    \draw (5*360/6:2.5)--({(4)*360/6+30}:3.5);

    \draw (4*360/6:2.5)--(5*360/6:2.5);
    
    \foreach \x in {0,2,4}{
      \draw[fill] (\x*360/6+30:3.5) circle (1.8pt);
      \draw (\x*360/6:3.5)--({\x*360/6}:2.5);
    }
    \foreach \x in {0,2,4}{
      \draw[fill] (\x*360/6:3.5) circle (1.8pt);

      \node at (1.20, 2.6) {$0$};
      \node at (2.95, 2.3) {$1$};
      \node at (3.45, 0.4) {$2$};
      \node at (2.85, 0.4) {$3$};
      \node at (1.20, -1.75) {$4$};
      \node at (0, -2.75) {$5$};
      \node at (-2.0, -2.75) {$6$};
       
      \node at (-1.10, -1.75) {$7$};
      \node at (-2.95, 0.4) {$8$};
      \node at (-3.3, 2.3) {$9$};

      \node at (-1.30, 3.3) {$10$};
      \node at (-1.10, 2.6) {$11$};
      
      \node at (0, 4.8) {$G''+e_5$};
      
      \draw[line width=1.5pt] (-0.2,0) circle [radius=4.3];
    }
  \end{scope}

  \begin{scope}[shift={(6,10)}]
    \foreach \x in {0,1,...,5}{
      \draw[fill] (\x*360/6:2.5) circle (1.8pt);
    }
    \draw (1*360/6:2.5)--({(1+1)*360/6}:2.5);
    \draw (2*360/6:2.5)--({(2+1)*360/6}:2.5);
    \draw (3*360/6:2.5)--({(3+1)*360/6}:2.5);
    \draw (5*360/6:2.5)--({(5+1)*360/6}:2.5);

    \draw (1*360/6:2.5)--({(0)*360/6+30}:3.5);
    \draw (0*360/6:2.5)--({(0)*360/6+30}:3.5);

    \draw (2*360/6:2.5)--({(2)*360/6+30}:3.5);
    \draw (3*360/6:2.5)--({(2)*360/6+30}:3.5);

    \draw (4*360/6:2.5)--({(4)*360/6+30}:3.5);
    \draw (5*360/6:2.5)--({(4)*360/6+30}:3.5);
    
    \foreach \x in {0,2,4}{
      \draw[fill] (\x*360/6+30:3.5) circle (1.8pt);
      \draw (\x*360/6:3.5)--({\x*360/6}:2.5);
    }
    \foreach \x in {0,2,4}{
      \draw[fill] (\x*360/6:3.5) circle (1.8pt);

      \node at (1.20, 2.6) {$0$};
      \node at (2.95, 2.3) {$1$};
      \node at (3.45, 0.4) {$2$};
      \node at (2.85, 0.4) {$3$};
      \node at (1.20, -1.75) {$4$};
      \node at (0, -2.75) {$5$};
      \node at (-2.0, -2.75) {$6$};
       
      \node at (-1.10, -1.75) {$7$};
      \node at (-2.95, 0.4) {$8$};
      \node at (-3.3, 2.3) {$9$};

      \node at (-1.30, 3.3) {$10$};
      \node at (-1.10, 2.6) {$11$};
      
      \node at (0, 4.8) {$G''$};
      
      \draw[line width=1.5pt] (-0.2,0) circle [radius=4.3];
    }
  \end{scope}

\end{tikzpicture}
        \caption{A part of the recursion tree showing the third situation that may lead to isomorphic copies.}\label{fig:situationThree}
\end{figure}
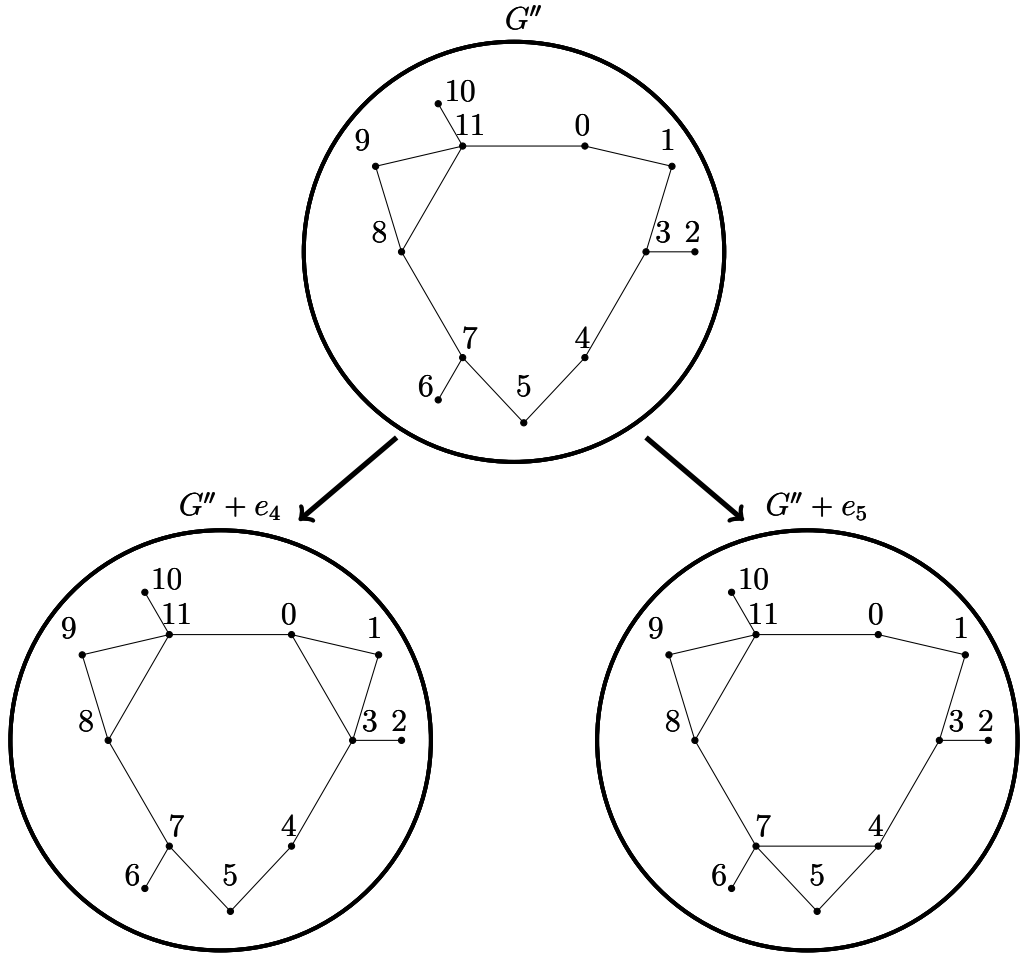
    
    The second and third source of isomorphism can be avoided by defining a so-called \textit{canonical reduction}. This is a reduction that is unique (up to isomorphism) and we will refer to the reverse operation as a \textit{canonical expansion}. The idea to avoid this source of isomorphism is to only expand a graph by doing canonical expansions. In other words, one only allows expansions such that their reverse operation is a canonical reduction. Note that if $G+e_1$ and $G'+e_3$ are isomorphic and $G$ and $G'$ are not isomorphic, then only one of the operations that removes $e_1$ or $e_3$ could have been a canonical reduction. Similarly, $G''+e_4$ and $G''+e_5$ are isomorphic, but only one of the operations that removes $e_4$ or $e_5$ could have been a canonical reduction. In practice, one can define such a canonical reduction by defining a tuple for each edge of a graph and choosing the edge that leads to the canonical reduction as the edge that is in the same edge orbit as the one that minimizes this tuple (in a lexicographical way). Here, each element of the tuple is the output of a function that can be evaluated for every edge of a graph. These functions are meant to be ordered from cheaply computable first (for example the sum of the degrees of the endpoints of the edge) to more expensive to compute last (for example the edge label in the canonical labeling of the graph). Often, this also corresponds to ordering the functions in increasing order of discriminating power. By doing this, one can find the edge with the lexicographically minimal tuple by computing the elements of the tuples one by one, which means that the more expensive functions often do not have to be evaluated if the first few elements of the tuple have a unique minimizing edge. Note that in all cases, it is necessary that the complete tuple has a unique minimizing edge (this is often achieved by having an element in the tuple that represents the edge label in the canonical labeling of the graph). The precise elements in the tuple are often problem-dependent and could exploit extra knowledge concerning the graph class that one is interested in generating.

    The \texttt{geng} algorithm from the \texttt{nauty} package~\cite{MP14} and \texttt{snarkhunter}~\cite{BGHM13,BGM11} are examples of graph generation algorithms that rely on the fourth technique.
    
    \item The fifth technique is generation based on the \textit{homomorphism principle}~\cite{GKL92}. This is a conceptually simple, yet powerful technique. Suppose we are interested in generating all graphs in a certain class. The main idea is to relate this class of graphs to a new class of graphs that is more coarse in such a way that any isomorphism between graphs in the first class also gives rise to an isomorphism between graphs in the second class. One then proceeds by first generating all graphs from the second (coarse) class and secondly one generates all graphs from the first class based on the graphs generated in the first step. Both of these generation steps can be treated independently and use for example one of the previously discussed techniques. For example, suppose we are interested in generating the class of graphs obtained by adding at most one pendant vertex to each vertex of a cubic graph. The homomorphism principle can be naturally applied to this problem: the coarser set of graphs could be taken as the underlying cubic graph (obtained by removing the pendant vertices). As we discussed in Section~\ref{sec:generationAlgorithms}, there exist several efficient generators for cubic graphs that one wants to generate during the first step. During the second step, one then has to add the pendant vertices using one's favorite graph generation technique.

    The algorithm described in~\cite{B25} is an example that uses the homomorphism principle.
\end{itemize}
\subsection{Further considerations}
\subsubsection{Efficiency}
Since these graph generation algorithms will be executed to generate graphs, the most important aspect is typically the actual running time that these algorithms consume rather than their asymptotic time complexity. This means that implementing these algorithms often involves quite a lot of algorithm engineering and sometimes asymptotically suboptimal algorithms are preferred. Typically, these graph generation algorithms are implemented in an efficient programming language such as \texttt{C} or \texttt{C++}. Moreover, a common practice is to make use of clever data structures such as bitsets (which can help to achieve constant factor improvements), SIMD instructions~\cite{F05}, appropriate compile flags that lead to executables optimized for efficiency and the use of code profilers to identify bottlenecks in the implementation of the algorithms.

Additionally, these graph generation algorithms are often executed at a large scale (e.g., in a high-performance computing environment) and therefore have to be parallelized. For instance, we refer to the papers~\cite{GJLSZ24},~\cite{AM24} and~\cite{BGM22} in which all computations took around 35, 80 and 125 CPU-years, respectively. For algorithms that use the canonical construction path method or orderly generation of graphs (or the homomorphism principle, depending on which generation technique is used in each of its steps), this can be quite easily achieved by splitting the computation into $m$ parts that should be processed by $m$ independent processes (note that most algorithms discussed in Section~\ref{sec:generationAlgorithms} indeed support this). This can be done by labeling all vertices in the recursion tree at a given recursion level $\ell$ from $0$ to $n-1$ (where $n$ is the number of vertices in the recursion tree at level $\ell$) and then let the process that should compute part $i$ only recurse from vertices in the recursion tree whose label is congruent to $i$ modulo $m$. Note that for these techniques the different processes do not need to communicate with each other for this idea to work, whereas parallelization is more cumbersome for algorithms that use other techniques to avoid generating isomorphic copies. We stress that for cases where there are a lot of graphs that have to be generated, the graphs are typically not written to memory (because of the large memory requirement), but instead they are written to the standard output and then immediately read from the standard input by a subsequent program that does some computation on these graphs (for example, checking whether there is a graph that yields a counterexample to a conjecture and writing this graph to the output if this is the case). This idea is very important to be able to use these algorithms at a large scale and is commonly known as~\textit{piping}.

In order to keep the size of the output small, the graphs are typically first converted to \texttt{graph6} format or \texttt{sparse6} format before writing them to the standard output. These formats encode a (specific labeling of a) graph as a short string by properly manipulating the adjacency matrix of a graph (we refer the interested reader to~\cite{NautyUserGuide} for more details on how these formats work precisely). For example, the Petersen graph can be succinctly represented as the following string in \texttt{graph6} format: \texttt{IsP@OkWHG}. The \texttt{graph6} format is intended to be used for all graphs, whereas the \texttt{sparse6} format is more efficient for sparse graphs. 

\subsubsection{Correctness}

We also would like to draw the reader's attention to the fact that one should take extra care for the correctness of the results when using a computer-assisted graph theory approach. The level of extra care that is required highly depends on the application. If the output of the algorithm is a graph for which the algorithm claims that it has some property $P$, it is often straightforward to manually verify this without the need for a computer (of course, depending on how intricate $P$ is). If on the other hand the result obtained by the algorithm consists of the claim that some graph does not exist, it is really important to take appropriate measures to ensure the correctness of the claim. First of all, it is always a good idea to formally prove the correctness of the algorithm (i.e., prove that the algorithm indeed generates all graphs from the given class of a given order). Secondly, one should pay extra attention to typical implementation mistakes such as uninitialized variables, memory leaks, overflows and variable shadowing using appropriate tools that can help to automatically detect such mistakes. For example, this can be achieved by compilers using appropriate compile flags and debugging tools such as \texttt{Valgrind}~\cite{NS07}. Thirdly, one should try to check for implementation mistakes by comparing the output of the algorithm with various known results from the literature that were manually obtained by a classical mathematical proof. Fourthly, one should aim at developing independent implementations of slower algorithms that potentially rely on different ideas, but are expected to produce the same output as the faster target algorithm. One should then try to use both algorithms to independently solve reasonably small problem instances and one should of course obtain the same output. Using such an approach, it is often possible to verify millions of cases and by being careful one can be very confident that the results are correct. Finally, if the application allows it, one can try to use a formally-verified theorem-proving system such as \texttt{Rocq}. However, we stress that for several applications where the computations take many CPU-years, it is likely unrealistic to re-run these computations using such a system due to performance limitations.

\section{Databases containing graphs}
\label{sec:databases}
In this section, we introduce four databases containing a large number of graphs (see also~\cite{mathbases} for an overview of mathematical databases). The key difference between the graph censuses that were discussed in Section~\ref{sec:generationAlgorithms} and databases containing graphs is that the latter also allow users to \textit{search} for graphs in the database that satisfy certain properties. To achieve this, these databases not only keep track of the graphs themselves, but also of a large list of graph invariants. We note that several of these graph invariants can be computed using software libraries such as \texttt{nauty}~\cite{MP14}, \texttt{SageMath}~\cite{sagemath}, \texttt{Mathematica}~\cite{mathematica} and \texttt{Maple}~\cite{maple}, whereas for other invariants it is necessary to write custom code to compute them.

Apart from the search functionality, most of the databases that we will discuss also support other useful features. These features make such databases very accessible, since researchers using the database do not have to write any code (in fact, we will later discuss in Section~\ref{sec:newResults} an example of how this can lead to new results). The databases are as follows:
\begin{itemize}
    \item The first database that we will discuss is \textit{DiscreteZOO}~\cite{BV20}. This is the most recent and conceptually simplest database of the four databases that will be discussed. It currently consists of all vertex-transitive graphs until order 31~\cite{HR20}, cubic vertex-transitive graphs until order 1280~\cite{PSV13} and cubic arc-transitive graphs until order 2048~\cite{CombinatorialDataConder}. For each of these graphs, the database stores 15 boolean invariants indicating whether a graph belongs to a given class or not (arc-transitive, bipartite, Cayley, distance-regular, distance-transitive, edge-transitive, Eulerian, Hamiltonian, M{\"o}bius ladder, overfull, partial cube, prism, split, split Praeger-Xu, strongly regular) and 7 numeric invariants (chromatic index, clique number, diameter, girth, order, number of connected components, number of triangles). The search functionality allows users to specify values that the corresponding invariants must have and the results can be sorted based on these invariants.

\item The second database is \textit{PHOEG}~\cite{DHM19}. It contains all pairwise non-isomorphic graphs on at most 10 vertices and stores for each of these graphs the value of 8 boolean invariants and 47 numeric invariants, many of which come from the domain of chemical graph theory such as the Wiener index~\cite{W47}, the Sombor index~\cite{G21} or the Randi{\'c} index~\cite{R75}. This database allows the users to filter graphs by specifying boolean values for the boolean invariants and lower and upper bounds for the numeric invariants. Subsequently, PHOEG plots these graphs as points in a two-dimensional space where the horizontal axis displays a first invariant and the vertical axis displays a second invariant (and each point can optionally receive a color given by a third invariant). The focus of PHOEG is mostly on extremal graph theory and therefore it also computes and displays the convex hull of this set of points, making it easier to identify graphs that are extremal for one of the invariants and also identify inequalities between the invariants (given by the facets of the polytope). Moreover, all graphs corresponding to the same point can also be simultaneously visualized (using different graph drawing algorithms) and an encoding of these graphs in \texttt{graph6} format is displayed. 

\item Thirdly, we discuss the \textit{Encyclopedia of Graphs}~\cite{EncyclopediaOfGraphs}. This database contains a very large variety of graph censuses (most of which have been mentioned in Section~\ref{sec:generationAlgorithms}) and also a large variety of graph invariants. The precise graph invariants that the database stores may be different for graphs belonging to different censuses (unlike the previous two databases). This also has implications for the search functionality: the database offers an option to search for graphs based on the graph's name, an internally stored Unique Graph Identifier (UGI) or a textual description of the collection in which the graph could be found. Within a fixed census, it is then possible to do a more refined search based on the graph invariants that are stored for the corresponding graph census. Additionally, this database allows users to add comments to graphs (in which one can for example refer to sources in which that particular graph plays an important role), download the graphs and their invariants and view an image of a selected graph.

\item Finally, the fourth database that we will discuss is the \textit{House of Graphs}~\cite{BCGM13,HOG}, which is to the best of our knowledge the most widely used database within the graph theory community.  One of the ideas behind this database is that some graphs repeatedly appear in various problems (for example as counterexamples to conjectures or as extremal graphs) and could therefore be considered to be more interesting than others. This database allows users to upload graphs that they find interesting (for example because they appear in a publication) and contains thousands of such graphs coming from various sources. Moreover, users can also add comments to these graphs that usually explain why a graph is interesting and point others to relevant sources, which is often useful to discover (sometimes unexpected) connections between different problems. For each graph, the database currently stores 12 boolean invariants (acyclic, bipartite, claw-free, connected, Eulerian, Hamiltonian, hypohamiltonian, hypotraceable, planar, regular, traceable, twin-free) and 39 numeric invariants, many of which require solving an \NP-hard problem. The search functionality that this database offers is the most advanced one out of the four databases. More precisely, users can form equalities or inequalities by combining multiple invariants that the graphs should adhere to. Therefore, this generalizes the search functionalities that the other databases offer with respect to the invariants. Additionally, users can also search for graphs associated with a given textual description (appearing for example in a comment), graphs associated with an internally stored House of Graphs identifier, graphs in which an input graph $H$ appears (or does not appear) as an (induced) subgraph, graphs with a given \texttt{graph6} representation or the $K$ most recently added graphs.

Search results can be sorted based on the invariants and once a graph is selected, the value of each graph invariant is shown. The graphs (and the values of their invariants) can also be downloaded in multiple formats. Moreover, the database also indicates whether the line graph or the complement of the current graph is already in the database or not. The database also displays multiple visualizations of the graph and users can add their own visualizations as well. We note that there are often multiple sensible ways to visualize a graph and seeing a graph visualized in the correct way can often tell one a lot about the structure one might be looking for.

Apart from these search-related features, the House of Graphs also allows users to draw new graphs (regardless of whether they are already in the database), modify existing visualizations and download them in several formats. Additionally, they also host and keep track of large censuses of graphs (not all of these are included in the searchable database to keep the computations feasible) and publications that interact with the House of Graphs.
\end{itemize}

\section{Other algorithmic paradigms}
\label{sec:otherParadigms}
As we have seen in the previous part, graph generation often plays an important role in a computer-assisted graph theory approach. Apart from graph generation, there are also many other algorithmic paradigms that complement graph generation. In this section, we discuss such paradigms and give some detailed examples from the literature.

\subsection{(Mixed integer) linear programming}
Linear programming is concerned with optimizing a linear objective function subject to linear constraints (equalities or inequalities) over a set of real decision variables. Without loss of generality, a linear programming problem with $n$ decision variables and $m$ constraints can be formulated as the following problem in standard form:
\[
\begin{aligned}
\text{Maximize} \quad & \bm{c}^\top \bm{x} \\
\text{subject to} \quad  A\bm{x} &\leq \bm{b} \\
                         \bm{x} &\geq \bm{0}
\end{aligned}
\]
Here, $\bm{x} \in \mathbb{R}^n$ is a vector containing the $n$ decision variables for which one should find the optimal value and the input for the problem consists of $\bm{c} \in \mathbb{R}^n$ (the vector of objective function coefficients),  $A \in \mathbb{R}^{m \times n}$ (the matrix of constraint coefficients) and $\bm{b} \in \mathbb{R}^m$ (the right-hand side vector of constraint bounds).

Using any computer-assisted graph theory approach requires one to take extra care that the programs do not contain any implementation mistakes since this could lead to incorrect outcomes. For linear programming problems, this problem can be easily solved based on the concept of \textit{duality}, which we will explain next. This allows one to present a certificate that the optimal solution was indeed found. In other words, this takes away the need for one to trust that the linear programming solver was implemented correctly. More precisely, one can associate a new linear programming problem related to the problem that is formulated in standard form above. This new linear programming problem looks as follows:
\[
\begin{aligned}
\text{Minimize} \quad & \bm{b}^\top \bm{y} \\
\text{subject to} \quad  A^\top \bm{y} &\geq \bm{c} \\
                         \bm{y} &\geq \bm{0}
\end{aligned}
\]
Here, $\bm{y} \in \mathbb{R}^m$ is a new vector containing $m$ decision variables. This new linear programming problem is referred to as the \textit{dual problem}, whereas the original problem is referred to as the \textit{primal problem}. We now introduce the weak duality theorem and strong duality theorem for linear programming.
\begin{thr}[Weak duality theorem for linear programming~\cite{BT97}]
If $\bm{x}$ is a feasible solution for the primal problem and $\bm{y}$ is a feasible solution for the dual problem, then $\bm{c}^\top \bm{x} \leq \bm{b}^\top \bm{y}$.
\end{thr}

\begin{thr}[Strong duality theorem for linear programming~\cite{BT97}]
If $\bm{x}^*$ is an optimal solution for the primal problem and $\bm{y}^*$ is an optimal solution for the dual problem, then $\bm{c}^\top \bm{x}^* = \bm{b}^\top \bm{y}^*$.
\end{thr}

This means that if one claims to have found an optimal solution $\bm{x}^*$ for the primal problem, one can prove this claim by presenting a feasible solution $\bm{y}$ for the dual problem such that $\bm{c}^\top \bm{x}^* = \bm{b}^\top \bm{y}$.

\textit{Mixed integer linear programming} (MILP) is a notion that is closely related to linear programming. However, in the case of MILP, one can also add constraints that indicate that a certain variable should belong to a certain set of integers. This is often very useful, for example to model that a certain vertex is chosen or not chosen or that a certain edge is present or not present using an indicator variable $(x \in \{0,1\})$. Linear programming problems can be solved in polynomial time using the ellipsoid method if all coefficients are rational numbers~\cite{K80}, but the most efficient practical algorithms are based on the simplex algorithm~\cite{D16} or interior-point methods~\cite{K84}. MILP on the other hand is a classical NP-hard problem~\cite{BT97}, but is often still tractable in practice. Modern (mixed integer) linear programming solvers such as \texttt{CPLEX}~\cite{CPLEX} and \texttt{Gurobi}~\cite{Gurobi} are often able to solve problems involving thousands of variables and constraints in a short amount of time to optimality. We refer the interested reader to~\cite{BT97,C83,Wolsey20} for a much deeper treatment of this rich area of research.

Several problems in graph theory can be modeled using MILP. We now discuss a recent example of how MILP was used by Wagner~\cite{Wagner20} for a problem in extremal graph theory. We first introduce some notations. For graphs $G$ and $H$, the Tur{\'a}n number $\operatorname{ex}(G,H)$ refers to the maximum number of edges that a subgraph of $G$ can have without containing $H$ as a subgraph. The \textit{complete $r$-partite graph} $K_{n_1,n_2,\ldots,n_r}$ is the graph consisting of $r$ independent sets $V_1,V_2,\ldots,V_r$ of sizes $n_1,n_2,\ldots,n_r$, respectively, such that there is an edge between each pair of vertices $u \in V_i$, $v \in V_j$ for which $i \neq j, i,j \in \{1,\ldots,r\}$. For an integer $k$ and a graph $G$, the graph $kG$ represents the disjoint union of $k$ copies of $G$. 

The authors in~\cite{DHKSY18} study $\operatorname{ex}(K_{n_1,n_2,\ldots,n_r},kK_s)$. One of their results is a construction from which they derive the following.
\begin{prop}[\cite{DHKSY18}]
For all integers $k, n_1, n_2, n_3, n_4$ such that $2 \leq k \leq n_1+n_2$, we have $\operatorname{ex}
(K_{n_1,n_2,n_3,n_4},kK_3) \geq (n_1+n_2+n_3)n_4+(k-1)n_3$.
\end{prop}
Wagner~\cite{Wagner20} used MILP to calculate $\operatorname{ex}(K_{n_1,n_2,n_3,n_4},kK_3)$ exactly for small parameters (together with the corresponding subgraph achieving this number of edges). For a given set of integers $n_1,n_2,n_3,n_4,k$, let $\mathcal{E}$ consist of all subsets of edges of $E(K_{n_1,n_2,n_3,n_4})$ that induce a copy of $kK_3$. One can now introduce an indicator variable $x_e \in \{0,1\}$ for each edge $e \in E(K_{n_1,n_2,n_3,n_4})$ and the value $\operatorname{ex}(K_{n_1,n_2,n_3,n_4},kK_3)$ is then given by the optimum of the following optimization problem:
\[
\begin{aligned}
\text{Maximize} \quad & \sum_{e \in E(K_{n_1,n_2,n_3,n_4})} x_e\\
\text{subject to} \quad  \sum_{e \in E'} x_e &\leq 3k-1, \quad E' \in \mathcal{E}\\
                         x_e &\in \{0,1\}, \quad e \in E(K_{n_1,n_2,n_3,n_4})
\end{aligned}
\]
By inspecting the optimal solution for small values and generalizing the pattern, Wagner~\cite{Wagner20} was able to improve the result in the case where $2 \leq k \leq n_1=n_2=n_3=n_4$.
\begin{prop}[\cite{Wagner20}]
For all integers $2 \leq k \leq n$, we have $\operatorname{ex}(K_{n,n,n,n},kK_3) \geq 4n^2+(k-1)n$.
\end{prop}

Apart from the previous example where MILP was used, we also mention that linear programming was recently used to automate a part of the \textit{discharging method}\footnote{Note that the proofs of the Four Color Theorem mentioned in the introduction of this survey also crucially rely on the discharging method, but no linear programming was used in those proofs.} (see~\cite{CW17} for a nice survey on the discharging method). The goal of using this method is to show that every graph $G$ from a certain class of graphs $\mathcal{G}$ has to satisfy some property $P$ (where typically $\mathcal{G}$ is a class of sparse graphs, for example planar graphs). 

More precisely, the method operates by first assigning \textit{charges} (real numbers) to elements of $G \in \mathcal{G}$ (typically the elements are vertices or faces) in such a way that the sum of all charges is negative. Subsequently, one defines a set of \textit{discharging rules} that allow one to redistribute the charges based on the properties of the elements (e.g., based on the vertex degrees) without affecting the sum of all charges. One then proceeds by defining a set of \textit{configurations} $\mathcal{C}$ and by concluding that each graph $G \in \mathcal{G}$ must \textit{contain} (the precise meaning of \textit{contain} depends on the specific use case) at least one of the configurations in $\mathcal{C}$, because assuming otherwise would result in a situation where the sum of all charges would be non-negative after applying the discharging rules. Finally, one considers a minimal graph $G \in \mathcal{G}$ (e.g., minimizing the pair $(|V(G)|,|E(G)|)$) that does not satisfy $P$ and shows that each configuration in $\mathcal{C}$ that appears in $G$ can be used to \textit{reduce} $G$ to a smaller graph $G' \in \mathcal{G}$ such that $G'$ also does not satisfy $P$ (such a configuration is called \textit{reducible}), thereby obtaining a contradiction with the assumption that there exists a (minimal) graph in $\mathcal{G}$ not satisfying $P$.

The authors of~\cite{BDDP24} recently used linear programming to conclude that either no set of suitable discharging rules exist to prove the desired property (under some assumptions about what the discharging rules can look like) or return a list of configurations such that at least one of them has to be reduced in order to obtain the contradiction discussed in the previous paragraph. Proving that configurations are reducible is something that can either be done manually or by using a heuristic (with the help of a computer). The authors of~\cite{BDDP24} iteratively use the linear programming approach (with new sets of constraints depending on the configurations) and prove new reducibility results until this iterative process stabilizes. The result from this iterative process is either a full discharging proof, or if one cannot reduce any output configuration, one has to rethink the modeling of the problem. By doing this, they made progress on Wegner's well-known conjecture~\cite{W77} about the distance 2-chromatic number $\chi_2(G)$ of planar graphs $G$ with maximum degree $\Delta$ (i.e., the classical chromatic number of the graph obtained by adding an edge between each pair of vertices at distance 2 in $G$):
\begin{conj}[\cite{W77}]
Let $G$ be a planar graph with maximum degree $\Delta$. Then, we have:
\begin{align*}
    \chi_2(G) \leq  \begin{cases}
        7& \textrm{if } \Delta=3,\\
        \Delta+5& \textrm{if } 4 \leq \Delta \leq 7,\\
        \lfloor \frac{3\Delta}{2} \rfloor +1 & \textrm{if } \Delta \geq 8.\\
    \end{cases}
\end{align*}
\end{conj}

More precisely, the authors of~\cite{BDDP24} improved the best known upper bound for the case $\Delta=4$:
\begin{thr}[\cite{BDDP24}]
Let $G$ be a planar graph with maximum degree $4$. Then, we have $\chi_2(G) \leq 12$.
\end{thr}
\subsection{Semidefinite programming and flag algebras}
\label{sec:semidefiniteProgrammingAndFlagAlgebras}
A symmetric matrix $M \in \mathbb{R}^{n \times n}$ is called \textit{positive semidefinite} if for all $\bm{x} \in \mathbb{R}^n$, we have $\bm{x}^\top M\bm{x} \geq 0$. Equivalently, a symmetric matrix $M \in \mathbb{R}^{n \times n}$ is called positive semidefinite if and only if all of its eigenvalues are nonnegative. This is denoted by $M \succeq 0$. Let $\mathbb{S}^{n}$ denote all symmetric matrices in $\mathbb{R}^{n \times n}$. Define the inner product of two symmetric matrices $A, B \in \mathbb{S}^{n}$ as $\langle A, B \rangle= \sum_{i,j=1}^{n}A_{i,j}B_{i,j}$. Now a semidefinite programming problem asks one to find a matrix $X \in \mathbb{S}^n$ that solves the following problem (given in standard form):
\[
\begin{aligned}
\text{Minimize} \quad & \langle C, X \rangle \\
\text{subject to} \quad & \langle A_i, X \rangle = b_i, \quad i \in \{1,\ldots,m\}\\
                \quad &X \succeq 0
\end{aligned}
\]
Here, the input variables are $C, A_1,\ldots,A_m \in \mathbb{S}^n$ and $\mathbf{b} = (b_1,\ldots,b_m) \in \mathbb{R}^m$.

Similarly to what we saw for linear programming problems, for semidefinite programming problems it is also possible to associate a new semidefinite programming problem. It is again called a \textit{dual problem}, whereas the original one is called a \textit{primal problem}. The dual problem of the primal problem mentioned above looks as follows (where $\mathbf{y} = (y_1,\ldots,y_m) \in \mathbb{R}^m$ should be determined):
\[
\begin{aligned}
\text{Maximize} \quad & \mathbf{b}^\top\mathbf{y} \\
\text{subject to} \quad &  C-\sum_{i=1}^{m} y_iA_i \succeq 0\\
                \quad &\mathbf{y} \in \mathbb{R}^m
\end{aligned}
\]

Again, we have a weak duality theorem for semidefinite programming.
\begin{thr}[Weak duality theorem for semidefinite programming~\cite{VB96}]
If $X$ is a feasible solution for the primal problem and $\bm{y}$ is a feasible solution for the dual problem, then $\langle C, X \rangle \geq \mathbf{b}^\top\mathbf{y}$.
\end{thr}

However, unlike the case of linear programming, for semidefinite programming problems we do not have an analogous \textit{strong duality theorem for semidefinite programming} (in general). In other words, there exist semidefinite programming problems such that if $X^*$ is an optimal solution for the primal problem and $\bm{y}^*$ is an optimal solution for the dual problem, then $\langle C, X^* \rangle > \mathbf{b}^\top\mathbf{y}^*$ (see~\cite{VB96} for more information on the subject). Semidefinite programming has several applications and efficient solvers such as \texttt{MOSEK}~\cite{MOSEK}, \texttt{SDPA}~\cite{SDPA} and the \texttt{CSDP} solver from the \texttt{COIN-OR} project~\cite{COIN-OR} exist.

One important application of semidefinite programming in the context of computer-assisted graph theory is \textit{flag algebras} introduced by Razborov~\cite{R07}. The flag algebra method is useful for proving bounds for problems in extremal graph theory that deal with subgraph densities. A \textit{flag} can be thought of as a (typically small) graph with labeled vertices, allowing one to fix a part of a structure. The method allows one to transform a combinatorial problem into algebraic inequalities involving densities of small flags. (We note that it is often necessary to exhaustively generate all small flags and the methods discussed in Section~\ref{sec:generationAlgorithms} are often helpful for this.) The main idea is then to encode these inequalities and the problem-specific constraints into a semidefinite programming problem whose optimum yields a bound on the density of a target subgraph. In this context we also wish to mention \texttt{Flagmatic}~\cite{Flagmatic}, which is a very useful software package that automates many steps of the flag algebra method. For a more detailed technical introduction to the theory of flag algebras, we refer the interested reader to~\cite{R07}. 

\begin{figure}[h]
    \centering

    \begin{tikzpicture}[scale=1.2]

  \def\radius{3}           
  \def\groupsep{0.6}       
  \def\numgroups{5}
  \def\pointspergroup{4}

  \foreach \i in {0,...,4} {
    \pgfmathsetmacro\angle{72*\i}
    \coordinate (center\i) at (\angle:\radius);
    
    \pgfmathsetmacro\dirangle{\angle + 90}
    
    \foreach \j in {0,...,3} {
      \pgfmathsetmacro\offset{(\j - 1.5) * \groupsep}
      \path let \p1 = (center\i) in
        coordinate (P\i\j) at ($(\p1) + (\dirangle:\offset)$);
      \fill (P\i\j) circle (1.5pt); 
    }
  }  
  
  \foreach \i in {0,...,4} {
    \pgfmathtruncatemacro\targeti{mod(\i+1, 5)}
    \foreach \j in {0,...,3} {
      \foreach \k in {0,...,3} {
        \draw (P\i\j) -- (P\targeti\k);
      }
    }
  }
      
\end{tikzpicture}
\caption{A triangle-free graph on $20$ vertices with exactly $(\frac{20}{5})^5=1024$ cycles of length $5$.}\label{fig:C5BlowUp}
\end{figure}
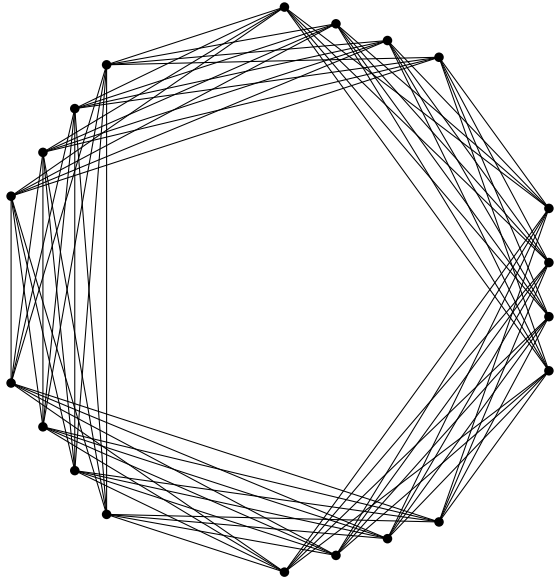

There are plenty of important problems in graph theory in which the flag algebra method was used, for instance see the papers~\cite{BHLPVY17,CKPSTY13,DHMNS13,G12,HHKNR13,LP18,LP21,R16}. We now discuss two of these problems. In 1984, Erd{\H{o}}s made the following conjecture.
\begin{conj}[\cite{E84}]
Any triangle-free graph on $n$ vertices contains at most $(\frac{n}{5})^5$ cycles of length $5$.
\end{conj}
The balanced blow-up of the $5$-cycle $C_5$ by independent sets (see Fig.~\ref{fig:C5BlowUp}) achieves this bound (when $n$ is divisible by $5$).
Erd{\H{o}}s also made the following conjecture.
\begin{conj}[\cite{E84}]
For graphs of order $n \geq 5$, the balanced blow-up of the $5$-cycle $C_5$ maximizes the number of 5-cycles in the class of triangle-free graphs on $n$ vertices (here, every vertex in $C_5$ is replaced by an independent set of size $\lfloor \frac{n}{5}\rfloor$ or $\lceil \frac{n}{5}\rceil$ and every edge is replaced by a complete bipartite graph between the corresponding sets).
\end{conj}

The authors of~\cite{G12} and~\cite{HHKNR13} obtained independent proofs of the claim that any triangle-free graph on $n$ vertices contains at most $(\frac{n}{5})^5$ cycles of length $5$. Let $C_{C_5}(G)$ denote the set of all $5$-element subsets of $V(G)$ inducing a copy of $C_5$, let $\operatorname{ex}_{C_5}(n,K_3)=\max\{|C_{C_5}(G)|: G\text{ is triangle-free and } |V(G)|=n\}$, let $\pi_{C_5}(K_3)= \lim_{n \rightarrow \infty}\frac{\operatorname{ex}_{C_5}(n,K_3)}{\binom{n}{5}}$ and let $\mathcal{H}$ denote the set of all triangle-free graphs on at most $5$ vertices. The crucial inequality derived from flag algebras is of the form $\pi_{C_5}(K_3) \leq f(C_5,\mathcal{H})$. Here, the right-hand side $f(C_5,\mathcal{H})$ is a function derived from the flag algebra method that can be bounded by solving a suitable semidefinite programming problem. Later, the authors of~\cite{LP18} showed that for all $n \geq 10$, the only triangle-free graphs on $n$ vertices maximizing the number of cycles of length $5$ are obtained by the balanced blow-up of $C_5$. The remaining cases for $n \leq 9$ were handled by exhaustive enumeration (and it turned out that the M{\"o}bius ladder on 8 vertices is another triangle-free graph maximizing the number of $5$-cycles).

The second problem involves \textit{Ramsey numbers}. Ramsey's Theorem~\cite{R30} implies that for any $k$ graphs $G_1,G_2,\ldots,G_k$, every edge-coloring using the colors from $\{1,\ldots,k\}$ of every large enough clique $K_n$ on $n$ vertices contains, for some $i \in \{1,\ldots,k\}$, a copy of the graph $G_i$ in which all edges receive color $i$. The Ramsey number $R(G_1,\ldots,G_k)$ is the smallest integer $n$ for which this property holds. The flag algebra method is designed for finding asymptotic results for very large graphs and therefore one might suspect at first glance that this method is not suitable for determining Ramsey numbers. However, it turns out that this intuition is not correct as shown in the paper by Lidick{\'y} and Pfender~\cite{LP21}. More precisely, they derive upper bounds on Ramsey numbers by providing lower bounds for the density of independent sets of size $\ell$ in the class of blow-ups of $k$-edge-colored cliques that do not contain a copy of the graph $G_i$ in which all edges receive color $i$ for all $i \in \{1,\ldots,k\}$. These density lower bounds can be found based on flag algebras (by solving a suitable semidefinite programming problem). By combining the obtained upper bounds on Ramsey numbers with lower bounds (mostly from existing literature), Lidick{\'y} and Pfender~\cite{LP21} showed the following equalities.
\begin{thr}[\cite{LP21}]
The following equalities hold: $R(K_4^-,K_4^-,K_4^-)=28$, $R(K_8,C_5)=29$, $R(K_9,C_6)=41$, $R(Q_3,Q_3)=13$, $R(K_{3,5},K_{1,6})=17$ and $R(C_3,C_5,C_5)= 17$. Here, $C_k$ denotes the cycle on $k$ vertices, $K_k$ the clique on $k$ vertices, $K_{k_1,k_2}$ the complete bipartite graph with partite sets of size $k_1$ and $k_2$, $K_k^-$ the graph obtained by removing a single edge from $K_k$ and $Q_k$ the $k$-dimensional hypercube.
\end{thr}
\subsection{Dynamic programming}
Dynamic programming~\cite{B54} is a very powerful algorithmic paradigm for solving problems that can be reduced to finding appropriate recursion relations. It decomposes a problem into various (sometimes overlapping) subproblems and the key idea is to compute the solutions to all subproblems at most once by storing them in a table and looking them up again when required without recalculating the values.

This paradigm also frequently comes back in the context of computer-assisted graph theory. As we have seen in the previous sections, it is often necessary to efficiently calculate certain graph invariants and dynamic programming is often useful to do this. A textbook example is calculating the size of a maximum independent set in a tree $T$. Using dynamic programming, this can be done in linear time. More precisely, for a tree $T$, fix any vertex $v$ that will serve as the root of the tree. For a vertex $u \in V(T)$, let $N_i(u) := \{w \in V(T)~|~d(u,w)=i\}$ be the set of all vertices at distance $i$ from $u$. For every vertex $u \in N_i(v)$, let $T_u$ be the tree rooted at vertex $u$, i.e., the graph induced by $\cup_{j=0}^{\infty} \left( N_j(u) \cap N_{i+j}(v) \right)$. Now for each vertex $u \in V(T)$ and each $j \in \{0,1\}$, define $f(u,j)$ as the size of a maximum independent set in $T_u$ such that vertex $u$ is not contained in this set if $j=0$ and such that vertex $u$ is contained in this set if $j=1$. If $u$ has degree at most one, then clearly $f(u,j)=j$. Otherwise, for vertex $u \in N_i(v)$, $f(u,j)$ can be recursively expressed as follows:

	\[
	f(u,j) =  \left\{
	\begin{array}{lr}
		\sum_{w \in N_1(u) \cap N_{i+1}(v)} \max(f(w,0), f(w,1)) , & j=0 \mbox{,}\\
		1+\sum_{w \in N_1(u) \cap N_{i+1}(v)} f(w,0), & j=1 \mbox{.}
	\end{array} 
	\right.
	\]
This is the case because if $u$ is included in an independent set, then by definition none of its neighbors can be included, whereas if $u$ is not included, then its neighbors either could or could not be included. Using this recursive formulation, it is possible to compute $f(u,j)$ ($u \in V(T)$, $j \in \{0,1\}$) in $O(|V(T)|)$ time using dynamic programming. By the definition of $f$, the size of a maximum independent set in $T$ is then given by $\max(f(v,0),f(v,1))$. We note that a maximum independent set with this size can then easily be obtained by following the recursion starting from $v$ and each time comparing which of $f(w,0)$ and $f(w,1)$ is larger.

Apart from computing graph invariants, dynamic programming can also be useful for the generation of graphs (or as a subroutine thereof, see e.g.~\cite{GR13}). We now discuss a recent example from a paper by Cambie and the author~\cite{CJ25} in which one of the problems that is studied concerns finding a graph $G$ of largest \textit{diameter} $d$ (i.e., $d=\text{diam}(G) := \max_{u,v \in V(G)} d(u,v)$) among all graphs of order $n$, minimum degree $\delta$ and chromatic number at most $\chi$. More precisely, they are interested in the asymptotic case (where the order $n$ goes to infinity) by determining lower and upper bounds for $c_{\delta,\chi}$: the smallest rational number such that $\text{diam}(G) \leq c_{\delta,\chi}n+O(1)$ for all graphs $G$ with order $n$, minimum degree $\delta$ and chromatic number at most $\chi$. It was shown in~\cite{CJ25} that determining $c_{\delta,\chi}$ can be reduced to finding an optimal gadget that will be used in a gluing-like operation to create graphs of large diameter (in comparison with their order) with minimum degree $\delta$ and chromatic number at most $\chi$. 

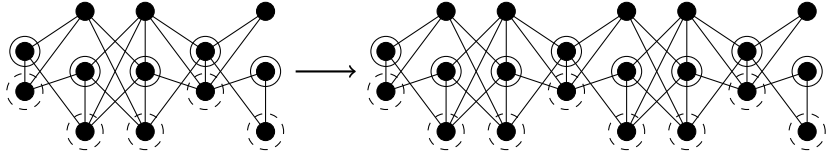
\begin{figure}[h]
    \centering

    \begin{tikzpicture}[scale=0.8]

    \foreach \x in {0,3}{
        \foreach \y in {-0.333,0.333}{
        \draw[fill] (\x,\y) circle (0.15);
        }
        \draw (\x,-0.333)--(\x,0.333);
    }
    \foreach \x in {2}{
        \foreach \y in {-1,0,1}{
        \draw[fill] (\x,\y) circle (0.15);
        }
        \draw (\x,-1)--(\x,0);
        \draw (\x,0)--(\x,1);
    }
    \foreach \x in {1,4}{
        \foreach \y in {-1,0,1}{
        \draw[fill] (\x,\y) circle (0.15);
        }
        \draw (\x,-1)--(\x,0);
    }
    \draw (1,-1)--(0,0.333);
    \draw (1,0)--(0,-0.333);
    \draw (1,1)--(0,0.333);
    \draw (1,1)--(0,-0.333);

    \draw (4,-1)--(3,0.333);
    \draw (4,0)--(3,-0.333);
    \draw (4,1)--(3,0.333);
    \draw (4,1)--(3,-0.333);

    \draw (2,1)--(1,0);
    \draw (2,1)--(1,-1);
    \draw (2,0)--(1,1);
    \draw (2,0)--(1,-1);
    \draw (2,-1)--(1,1);
    \draw (2,-1)--(1,0);

    \draw (3,0.333)--(2,-1);
    \draw (3,0.333)--(2,1);
    \draw (3,-0.333)--(2,0);
    \draw (3,-0.333)--(2,1);
   \draw[->,thick] (4.5,0) -- ++(1,0);

    \foreach \x in {6,9}{
        \foreach \y in {-0.333,0.333}{
        \draw[fill] (\x,\y) circle (0.15);
        }
        \draw (\x,-0.333)--(\x,0.333);
    }
    \foreach \x in {8}{
        \foreach \y in {-1,0,1}{
        \draw[fill] (\x,\y) circle (0.15);
        }
        \draw (\x,-1)--(\x,0);
        \draw (\x,0)--(\x,1);
    }
    \foreach \x in {7,10}{
        \foreach \y in {-1,0,1}{
        \draw[fill] (\x,\y) circle (0.15);
        }
        \draw (\x,-1)--(\x,0);
    }

    \draw (7,-1)--(6,0.333);
    \draw (7,0)--(6,-0.333);
    \draw (7,1)--(6,0.333);
    \draw (7,1)--(6,-0.333);

    \draw (10,-1)--(9,0.333);
    \draw (10,0)--(9,-0.333);
    \draw (10,1)--(9,0.333);
    \draw (10,1)--(9,-0.333);

    \draw (8,1)--(7,0);
    \draw (8,1)--(7,-1);
    \draw (8,0)--(7,1);
    \draw (8,0)--(7,-1);
    \draw (8,-1)--(7,1);
    \draw (8,-1)--(7,0);

    \draw (9,0.333)--(8,-1);
    \draw (9,0.333)--(8,1);
    \draw (9,-0.333)--(8,0);
    \draw (9,-0.333)--(8,1);

    \foreach \x in {12}{
        \foreach \y in {-0.333,0.333}{
        \draw[fill] (\x,\y) circle (0.15);
        }
        \draw (\x,-0.333)--(\x,0.333);
    }
    \foreach \x in {11}{
        \foreach \y in {-1,0,1}{
        \draw[fill] (\x,\y) circle (0.15);
        }
        \draw (\x,-1)--(\x,0);
        \draw (\x,0)--(\x,1);
    }
    \foreach \x in {13}{
        \foreach \y in {-1,0,1}{
        \draw[fill] (\x,\y) circle (0.15);
        }
        \draw (\x,-1)--(\x,0);
    }

    \draw (13,-1)--(12,0.333);
    \draw (13,0)--(12,-0.333);
    \draw (13,1)--(12,0.333);
    \draw (13,1)--(12,-0.333);

    \draw (11,1)--(10,0);
    \draw (11,1)--(10,-1);
    \draw (11,0)--(10,1);
    \draw (11,0)--(10,-1);
    \draw (11,-1)--(10,1);
    \draw (11,-1)--(10,0);

    \draw (12,0.333)--(11,-1);
    \draw (12,0.333)--(11,1);
    \draw (12,-0.333)--(11,0);
    \draw (12,-0.333)--(11,1);

    \draw[] (0,0.333) circle [radius=0.25];
    \draw[] (3,0.333) circle [radius=0.25];
    \draw[] (6,0.333) circle [radius=0.25];
    \draw[] (9,0.333) circle [radius=0.25];
    \draw[] (12,0.333) circle [radius=0.25];
    
    \draw[dashed] (0,-0.333) circle [radius=0.30];
    \draw[dashed] (3,-0.333) circle [radius=0.30];
    \draw[dashed] (6,-0.333) circle [radius=0.30];
    \draw[dashed] (9,-0.333) circle [radius=0.30];
    \draw[dashed] (12,-0.333) circle [radius=0.30];
    
    \draw[] (1,0) circle [radius=0.25];
    \draw[] (2,0) circle [radius=0.25];
    \draw[] (4,0) circle [radius=0.25];
    \draw[] (7,0) circle [radius=0.25];
    \draw[] (8,0) circle [radius=0.25];
    \draw[] (10,0) circle [radius=0.25];
    \draw[] (11,0) circle [radius=0.25];
    \draw[] (13,0) circle [radius=0.25];

    \draw[dashed] (1,-1) circle [radius=0.30];
    \draw[dashed] (2,-1) circle [radius=0.30];
    \draw[dashed] (4,-1) circle [radius=0.30];
    \draw[dashed] (7,-1) circle [radius=0.30];
    \draw[dashed] (8,-1) circle [radius=0.30];
    \draw[dashed] (10,-1) circle [radius=0.30];
    \draw[dashed] (11,-1) circle [radius=0.30];
    \draw[dashed] (13,-1) circle [radius=0.30];
    
    \end{tikzpicture} 
\caption{A repeatable graph with respect to $\delta=4$ and $\chi=3$ and a gluing-like operation. The different vertex styles show a proper $3$-coloring.}\label{fig:repeatableGraph}
\end{figure}

This gadget is called a \textit{repeatable graph} and its definition is chosen in such a way that this gluing-like operation has these properties. We recall the necessary definitions from~\cite{CJ25}: consider a graph $G$ and a subset of vertices $N_0 \subset V(G)$. For all $i \geq 1$, define the neighborhoods $N_i := \{v \in V(G)~|~\min_{u \in N_0}(d(u,v))=i\}$ and let $d$ be the largest integer for which $N_d$ is non-empty. Let $G[S]$ be the graph induced by the vertices in $S$. For integers $\delta$ and $\chi$, we call $G$ \textit{repeatable} with respect to $\delta$ and $\chi$ if $d \geq 2$ and $G[N_0 \cup N_1] \cong G[N_{d-1} \cup N_d]$ such that the isomorphism maps vertices in $N_0$ to $N_{d-1}$ and vertices in $N_1$ to $N_d$ and $G$ has chromatic number at most $\chi$ and every vertex in $G$, except for vertices in $N_0 \cup N_d$, has minimum degree at least $\delta$. Fig.~\ref{fig:repeatableGraph} shows an example of a repeatable graph with respect to $\delta=4$ and $\chi=3$ and the gluing-like operation that was mentioned before.

It was shown in~\cite{CJ25} that if $G$ is a repeatable graph with respect to $\delta$ and $\chi$ induced by neighborhoods $N_0, N_1, \ldots, N_d$, then $c_{\delta,\chi} \geq \frac{d-1}{|N_0 \cup N_1 \cup \ldots \cup N_{d-2}|}$ (and in fact it is even guaranteed that a repeatable graph must exist for which this inequality becomes an equality). Hence, it is of interest to generate repeatable graphs for which this ratio is large. Consider a proper $\chi$-coloring $c$ of a repeatable graph $G$ with respect to $\delta$ and $\chi$. The following heuristics make this search process easier (for some of these heuristics one can in fact prove that it must hold for every optimal repeatable graph):
\begin{itemize}
    \item If $u, v \in N_{i-1} \cup N_i$ and $c(u) \neq c(v)$, then one may assume that $uv \in E(G)$,
    \item $d \leq c'$ for some (small) constant $c'$,
    \item for all $i \in \{0,\ldots,d\}$, $|N_i| \leq c''$ for some (small) constant $c''$.
\end{itemize}
Bearing these heuristics in mind, finding a good repeatable graph naturally lends itself to a dynamic programming approach where one simultaneously searches for a good repeatable graph and a proper $\chi$-coloring $c$ thereof. More precisely, one can define a function $f(s,C_0,C_1,C_{s-1},C_{s})$ to indicate the least order $n$ such that there exists a graph induced by the neighborhoods $N_0, N_1, \ldots, N_{s-1}, N_s$ (where $|N_0 \cup \ldots \cup N_s|=n$) such that the number of colors in each color class of $c$ in $N_i$ is given by $C_i$ ($i \in \{0,1,s-1,s\}$) and every vertex in $N_1 \cup \ldots \cup N_{s-1}$ has degree at least $\delta$. Now $f(s,C_0,C_1,C_{s-1},C_{s})$ can be recursively expressed in terms of $f(s-1,C_0,C_1,C_{s-2},C_{s-1})$ by iterating over all possibilities for $C_{s-2}$ containing at most $c'$ vertices and adding the appropriate edges as determined by the $\chi$-coloring (the details are described in~\cite{CJ25}). By implementing this dynamic programming approach, the following result (among other things) was shown in~\cite{CJ25}.
\begin{prop}[\cite{CJ25}]
The following inequality holds: $c_{3,16} \geq \frac{31}{216}$.
\end{prop}
We note that the corresponding repeatable graph has $219$ vertices. Further arguments discussed in~\cite{CJ25} imply that this yields a counterexample to the following conjecture by Erd{\H{o}}s, Pach, Pollack and Tuza~\cite{EPPT89} in which the upper bound on the chromatic number is replaced by an upper bound on the clique number (in a regime where the conjecture was still open\footnote{The authors of~\cite{CSS21} already gave counterexamples to this conjecture, but left the regime $(r-1)(3r+2) \leq \delta \leq 2(r-1)(3r+2)(2r-3)$ still open.}):
\begin{conj}[\cite{EPPT89}]
\label{conj:EPPT_wrong}
Let $r, \delta \geq 2$ be fixed integers and let $G$ be a connected graph of order $n$ and minimum degree $\delta$.

\begin{itemize}
    \item[(i)] If $G$ is $K_{2r}$-free and $\delta$ is a multiple of $(r - 1)(3r + 2)$, then, as $n \to \infty$,
    \[
    \text{diam}(G) \leq \frac{2(r - 1)(3r + 2)}{(2r^2 - 1)} \cdot \frac{n}{\delta} + O(1)
    \]
    \[
    = \left( 3 - \frac{2}{2r - 1} - \frac{1}{(2r - 1)(2r^2 - 1)} \right) \frac{n}{\delta} + O(1).
    \]
    
    \item[(ii)] If $G$ is $K_{2r+1}$-free and $\delta$ is a multiple of $3r - 1$, then, as $n \to \infty$,
    \[
    \text{diam}(G) \leq \frac{3r - 1}{r} \cdot \frac{n}{\delta} + O(1)
    \]
    \[
    = \left( 3 - \frac{2}{2r} \right) \frac{n}{\delta} + O(1).
    \]
\end{itemize}
\end{conj}
\subsection{SAT-based approaches}
The satisfiability problem (abbreviated as SAT) asks whether a given Boolean formula can be satisfied by assigning a truth value (true or false) to each variable. It is the first problem for which it was proven that it is \NP-complete~\cite{K72}. Many practical and theoretical problems can be reformulated as SAT problems and this resulted in tremendous research efforts and very effective SAT solvers such as \texttt{CaDiCaL}~\cite{B18} and \texttt{Kissat}~\cite{B20}. If a Boolean formula can be satisfied, this can be easily shown by giving a satisfying truth assignment. However, many modern SAT solvers also support \textit{proof logging} (for example in the \texttt{DRAT} format or a trimmed version thereof based on \texttt{DRAT-trim}~\cite{WHH14}). This serves as a proof that a Boolean formula is not satisfiable and can be checked by formally-verified theorem-proving systems such as \texttt{ACL2}, \texttt{Rocq} and \texttt{Isabelle}~\cite{CHHKS17,K99}. 

We now give two examples of how SAT-based approaches were used in the context of computer-assisted graph theory. The first example is related to Keller's conjecture~\cite{K30}, which at first glance is not immediately related to graph theory. The conjecture states that if the $n$-dimensional Euclidean space is tiled by identical hypercubes, then there are two hypercubes that have an $(n-1)$-dimensional face in common. The conjecture was known to be true for $n \leq 6$ since 1940~\cite{P40a,P40b} and false for $n \geq 8$ since 2002~\cite{M02}, leaving only the case $n=7$ open. The relationship with graph theory is that the conjecture is false for $n=7$ if and only if at least one of the graphs $G_{7,3}, G_{7,4}$ or $G_{7,6}$ contains a clique of size $128$~\cite{BHMN22}. Here, $G_{n,s}$ is a Keller graph~\cite{CS90} having vertices $\{0,1,\ldots,2s-1\}^n$  and there is an edge between two vertices if and only if they differ in at least two coordinates and they differ by exactly $s$ in at least one coordinate. The authors of~\cite{BHMN22} showed that none of the three graphs have a clique of size $128$ by reducing this problem to a SAT problem and using symmetry-breaking techniques to speed up the search (all the proofs of unsatisfiability combined ended up having a size of $224$ gigabytes in the binary \texttt{DRAT} format).

The second example is a SAT-based approach for generating all graphs (up to isomorphism) within certain graph classes. Kirchweger and Szeider~\cite{KS24} recently proposed the SAT Modulo Symmetries framework. An important concept in this framework applied to graph generation is the concept of a \textit{partial graph} consisting of a set $V$ of vertices and a symmetric function $f(u,v)$ mapping every pair of distinct vertices $u \neq v \in V$ to an element in $\{0,1,2\}$. This function can be thought of as an edge not being present ($f(u,v)=0$), an edge being present ($f(u,v)=1$) or an edge being undecided ($f(u,v)=2$). Hence, graphs are partial graphs in which every pair of distinct vertices is mapped to either $0$ or $1$. The property describing the class of graphs to be generated should be encoded as a Boolean formula. The framework then generates the partial graphs using a SAT solver and a symmetry propagator that informs the SAT solver about the lexicographical minimality of the partially generated graphs. In case a partially generated graph is not lexicographically minimal, the Boolean formula for which the SAT solver is trying to find a model is altered to change the search space. Two of the main strengths of this framework are that it leads to a highly extensible graph generator and the framework can also benefit from the proof logging capabilities of the SAT solvers.
\subsection{Metaheuristics}

Metaheuristics~\cite{G86} offer a set of guidelines to develop heuristic algorithms for combinatorial optimization problems. While metaheuristics typically do not offer any theoretical guarantees, they are the leading paradigm for finding good (not necessarily optimal) solutions for such problems in a short amount of time. Metaheuristics perform a series of operations on solutions in order to find better ones. The nature of these operations gives rise to several different classes of metaheuristics such as constructive metaheuristics, local search-based metaheuristics, population-based metaheuristics and hybrid metaheuristics (see for example the reference work~\cite{GK03}).

We now explain the widely used metaheuristic Variable Neighborhood Search~\cite{MH97} based on the description in~\cite{CH00}. This metaheuristic works with two sets of neighborhood structures $N$ and $N'$, where $N$ will be used in the outer algorithm and $N'$ will be used in the inner algorithm (which is called from the outer algorithm). Here, the neighborhoods in $N$ are typically large (making exhaustive exploration difficult) and also nested (e.g., the first neighborhood $N_1$ contains all graphs obtained from removing exactly one edge from a graph $G$, the second neighborhood $N_2$ contains all graphs obtained from removing exactly two edges from a graph $G$ and so on). The neighborhoods in $N'$ on the other hand are typically small enough to make complete enumeration feasible. The inner algorithm (shown in Algorithm~\ref{algo:VNLS}) tries to improve an incumbent solution by consecutively looking at the best neighbor of the incumbent solution in the neighborhood $N'_k$ (for increasing $k$) until a local optimum is reached and no improvement can be found. The outer algorithm (shown in Algorithm~\ref{algo:VNS}) tries to improve an incumbent solution by consecutively looking at a random neighbor in neighborhood $N_k$ (for increasing $k$). It then calls the inner algorithm starting from this random neighbor. If the new solution improves the incumbent solution, it is replaced and the search is restarted. This process is repeated until a stopping criterion is met (e.g., the time limit is reached). 

\begin{algorithm}[ht!]
\caption{variableNeighborhoodLocalSearch(Solution $x'$, Integer $k'_{\text{max}}$, Set of $k'_{\text{max}}$ neighborhood structures $N'$)}
\label{algo:VNLS}
  \begin{algorithmic}[1]
        \STATE $k \gets 1$
        \STATE $\text{improved} \gets \FALSE$
        \WHILE{$k \leq k'_{\text{max}}$}
            \STATE $x'' \gets \text{bestSolutionFrom}(N'_k(x'))$
            \IF{$x''$ is better than $x'$}
                \STATE $x' \gets x''$
                \STATE $\text{improved} \gets \TRUE$
            \ENDIF
            \STATE $k \gets k+1$
            \IF{$k>k'_{\text{max}}$ and $\text{improved}=\TRUE$}
                \STATE $k \gets 1$
                \STATE $\text{improved} \gets \FALSE$
            \ENDIF
        \ENDWHILE
        \RETURN $x'$
  \end{algorithmic}
\end{algorithm}

\begin{algorithm}[ht!]
\caption{variableNeighborhoodSearch(Integer $k_{\text{max}}$, Set of $k_{\text{max}}$ neighborhood structures $N$, Integer $k'_{\text{max}}$, Set of $k'_{\text{max}}$ neighborhood structures $N'$)}
\label{algo:VNS}
  \begin{algorithmic}[1]
        \STATE $x \gets \text{getInitialSolution()}$
        \WHILE{Stopping criterion not met}
            \STATE $k \gets 1$
            \WHILE{$k \leq k_{\text{max}}$}
                \STATE $x' \gets \text{randomSolutionFrom}(N_k(x))$
                \STATE $x'' \gets \text{variableNeighborhoodLocalSearch}(x',k'_{\text{max}},N')$
                \IF{$x''$ is better than $x$}
                    \STATE $x \gets x''$
                    \STATE $k \gets 1$
                \ELSE
                    \STATE $k \gets k+1$
                \ENDIF
            \ENDWHILE
        \ENDWHILE
        \RETURN $x$
  \end{algorithmic}
\end{algorithm}

A prominent example in computer-assisted graph theory that uses a metaheuristic (more precisely, the metaheuristic Variable Neighborhood Search) is the \texttt{AutoGraphiX} system~\cite{CH00}. This system is mainly used in extremal graph theory applications where one tries to find a graph from a given class that minimizes or maximizes a graph invariant. This is done by posing the problem as an unconstrained combinatorial optimization problem (all constraints are moved to the objective function to be optimized by introducing suitable penalty functions for the constraints). The system then tries to find an extremal graph by applying the Variable Neighborhood Search metaheuristic to this unconstrained combinatorial optimization problem. It implements a few standard graph operations as neighborhood structures such as removing an edge, adding an edge and swapping an edge and a non-edge. Apart from finding extremal graphs, the system can also be used to find counterexamples to conjectures of the form $A(G) \leq B(G)$ (where $A(G)$ and $B(G)$ are graph invariants). It can do this by trying to find a graph that maximizes $A(G)-B(G)$ or it can be used to formulate new conjectures of this form. This system was later also updated to the \texttt{AutoGraphiX2} system~\cite{ABFCHHLM06} in which several modifications were made (e.g., by using a different local search strategy in the inner algorithm) and to the \texttt{AutoGraphiX3} system~\cite{C17} (which, for example, supports multi-objective optimization). Among others, these systems have been used to disprove conjectures of the automated conjecture making system \texttt{Graffiti}~\cite{F88,HC00} and study extremal graph theory problems involving the largest eigenvalue of the adjacency matrix of a graph~\cite{ABCHRSS08}, the algebraic connectivity~\cite{BDHO05} and the irregularity~\cite{HM05}.

Another example in which metaheuristics were effectively used is the paper by Exoo~\cite{E12} in which he used the metaheuristics Simulated Annealing~\cite{KGV83} and Tabu Search~\cite{G89} to prove the following bound on the classical two-color Ramsey number $R(4,6)$:
\begin{prop}[\cite{E12}]
    The following inequality holds: $R(4,6) \geq 36$.
\end{prop}

Using the notation of Section~\ref{sec:semidefiniteProgrammingAndFlagAlgebras}, we have $R(K_4,K_6) \geq 36$. This is proven by finding a graph on $35$ vertices that does not contain $K_4$ as a subgraph nor an independent set of size $6$. A straightforward choice to model this as a combinatorial optimization problem would be to try to find a graph on $35$ vertices that minimizes the number of subgraphs isomorphic to $K_4$ plus the number of subgraphs isomorphic to $6K_1$. However, this paper illustrates that sometimes it is worth using a \textit{surrogate function} as objective function: a function that does not directly measure the objective that one is interested in, but rather a slight variation, thereby driving the search away from known parts of the search space. For the $R(4,6)$ case, this was achieved in one of the approaches discussed in~\cite{E12} by adding a term to the objective function such that one tries to maximize the number of induced copies of the four-vertex path $P_4$.

\subsection{Approaches based on machine learning}
Given that machine learning has been widely adopted in almost all scientific disciplines, it is not surprising to see that it also has applications in computer-assisted graph theory. One popular paradigm in machine learning is \textit{reinforcement learning}~\cite{SB98}. In this paradigm, an algorithm is able to interact with an environment by taking \textit{actions} from a given set called the \textit{action space}. Each time an action is taken, the state of the environment changes and after taking sufficiently many actions, a reward is received. The goal for the algorithm is to learn which actions should be taken given the state of the environment in order to maximize the cumulative reward. Hence, the performance of the algorithm improves over time.

An example of a reinforcement learning algorithm is the \textit{deep cross-entropy method}. This algorithm was used in~\cite{W21} to find counterexamples to several open conjectures in graph theory. The deep cross-entropy method uses a neural network that receives as input a state and outputs a probability distribution representing the likelihood of each possible action to be taken in the current state. Here, a graph on $n$ vertices can be encoded as a string consisting of $\frac{n(n-1)}{2}$ characters that can be $0$ or $1$. Each character corresponds to a pair of distinct vertices for which an edge is not present ($0$) or present ($1$). In the case of~\cite{W21}, the goal is to predict each of the $\frac{n(n-1)}{2}$ characters one by one such that the graph corresponding to the final encoding represents a counterexample to a conjecture. This is achieved by giving as input to the neural network a state that represents the first $k-1$ characters that were already predicted and a vector representing that the neural network currently has to predict character $k$. The output should represent a probability distribution that character $k$ should become a $0$ or $1$ (i.e., an action represents the decision to be made for character $k$). After all $\frac{n(n-1)}{2}$ characters have been predicted by sampling from the probability distributions, the neural network receives feedback by means of a reward function. This reward function should reflect how close the graph from the output is to being a counterexample. The algorithm operates by repeatedly generating $N$ graphs using the current weights of the neural network and then selecting the top $x \%$ of these graphs (where $N$ and $x$ are parameters) and updating the weights of the neural network in such a way that these graphs are more likely to be generated in the future. Therefore, the algorithm gradually improves over time and can effectively learn to generate graphs that are (close to being) counterexamples. The interested reader is referred to~\cite{L20} for a deeper treatment of neural networks and the deep cross-entropy method.

\begin{figure}[h]
    \centering

    \begin{tikzpicture}[scale=0.8]

     \draw[fill] (0,0) circle (0.15);
     \draw[fill] (3,0) circle (0.15);
     \draw[fill] (6,0) circle (0.15);
     \draw[fill] (9,0) circle (0.15);
     \draw[fill] (12,0) circle (0.15);
     \draw[fill] (15,0) circle (0.15);
     \draw[fill] (18,0) circle (0.15);
     \draw (0,0)--(3,0);
     \draw (3,0)--(6,0);
     \draw (6,0)--(9,0);
     \draw (9,0)--(12,0);
     \draw (12,0)--(15,0);
     \draw (15,0)--(18,0);

    
    \draw ($(6,0)+(26:2.3)$) -- ($(6,0)+(46:2.3)$);
    \draw ($(6,0)+(62:2.3)$) -- ($(6,0)+(82:2.3)$);
    \draw ($(6,0)+(98:2.3)$) -- ($(6,0)+(118:2.3)$);
    \draw ($(6,0)+(134:2.3)$) -- ($(6,0)+(154:2.3)$);


    \draw ($(6,0)$) -- ($(6,0)+(46:2.3)$);
    \draw ($(6,0)$) -- ($(6,0)+(82:2.3)$);
    \draw ($(6,0)$) -- ($(6,0)+(118:2.3)$);
    \draw ($(6,0)$) -- ($(6,0)+(154:2.3)$);

    \draw ($(6,0)+(26:2.3)$) -- ($(6,0)$);
    \draw ($(6,0)+(62:2.3)$) -- ($(6,0)$);
    \draw ($(6,0)+(98:2.3)$) -- ($(6,0)$);
    \draw ($(6,0)+(134:2.3)$) -- ($(6,0)$);
    
    \draw[fill] (6,0) ++(26:2.3) circle (0.15);
    \draw[fill] (6,0) ++(46:2.3) circle (0.15);
    
    \draw[fill] (6,0) ++(62:2.3) circle (0.15);
    \draw[fill] (6,0) ++(82:2.3) circle (0.15);

    \draw[fill] (6,0) ++(98:2.3) circle (0.15);
    \draw[fill] (6,0) ++(118:2.3) circle (0.15);

    \draw[fill] (6,0) ++(134:2.3) circle (0.15);
    \draw[fill] (6,0) ++(154:2.3) circle (0.15);


    \draw ($(6,0)+(180+20:2.7)$) -- ($(6,0)+(180+40:2.7)$);
    \draw ($(6,0)+(180+50:2.7)$) -- ($(6,0)+(180+70:2.7)$);
    \draw ($(6,0)+(180+80:2.7)$) -- ($(6,0)+(180+100:2.7)$);
    \draw ($(6,0)+(180+110:2.7)$) -- ($(6,0)+(180+130:2.7)$);
    \draw ($(6,0)+(180+140:2.7)$) -- ($(6,0)+(180+160:2.7)$);

    \draw ($(6,0)+(180+30:1.8)$) -- ($(6,0)+(180+40:2.7)$);
    \draw ($(6,0)+(180+60:1.8)$) -- ($(6,0)+(180+70:2.7)$);
    \draw ($(6,0)+(180+90:1.8)$) -- ($(6,0)+(180+100:2.7)$);
    \draw ($(6,0)+(180+120:1.8)$) -- ($(6,0)+(180+130:2.7)$);
    \draw ($(6,0)+(180+150:1.8)$) -- ($(6,0)+(180+160:2.7)$);

    \draw ($(6,0)+(180+20:2.7)$) -- ($(6,0)+(180+30:1.8)$);
    \draw ($(6,0)+(180+50:2.7)$) -- ($(6,0)+(180+60:1.8)$);
    \draw ($(6,0)+(180+80:2.7)$) -- ($(6,0)+(180+90:1.8)$);
    \draw ($(6,0)+(180+110:2.7)$) -- ($(6,0)+(180+120:1.8)$);
    \draw ($(6,0)+(180+140:2.7)$) -- ($(6,0)+(180+150:1.8)$);

    
    \draw ($(6,0)+(180+20:2.7)$) -- ($(6,0)$);
    \draw ($(6,0)+(180+50:2.7)$) -- ($(6,0)$);
    \draw ($(6,0)+(180+80:2.7)$) -- ($(6,0)$);
    \draw ($(6,0)+(180+110:2.7)$) -- ($(6,0)$);
    \draw ($(6,0)+(180+140:2.7)$) -- ($(6,0)$);

    \draw ($(6,0)+(180+30:1.8)$) -- ($(6,0)$);
    \draw ($(6,0)+(180+60:1.8)$) -- ($(6,0)$);
    \draw ($(6,0)+(180+90:1.8)$) -- ($(6,0)$);
    \draw ($(6,0)+(180+120:1.8)$) -- ($(6,0)$);
    \draw ($(6,0)+(180+150:1.8)$) -- ($(6,0)$);

    \draw ($(6,0)+(180+40:2.7)$) -- ($(6,0)$);
    \draw ($(6,0)+(180+70:2.7)$) -- ($(6,0)$);
    \draw ($(6,0)+(180+100:2.7)$) -- ($(6,0)$);
    \draw ($(6,0)+(180+130:2.7)$) -- ($(6,0)$);
    \draw ($(6,0)+(180+160:2.7)$) -- ($(6,0)$);
    
    \draw[fill] (6,0) ++(180+20:2.7) circle (0.15);
    \draw[fill] (6,0) ++(180+40:2.7) circle (0.15);
    
    \draw[fill] (6,0) ++(180+50:2.7) circle (0.15);
    \draw[fill] (6,0) ++(180+70:2.7) circle (0.15);

    \draw[fill] (6,0) ++(180+80:2.7) circle (0.15);
    \draw[fill] (6,0) ++(180+100:2.7) circle (0.15);

    \draw[fill] (6,0) ++(180+110:2.7) circle (0.15);
    \draw[fill] (6,0) ++(180+130:2.7) circle (0.15);

    \draw[fill] (6,0) ++(180+140:2.7) circle (0.15);
    \draw[fill] (6,0) ++(180+160:2.7) circle (0.15);

    \draw[fill] (6,0) ++(180+30:1.8) circle (0.15);
    \draw[fill] (6,0) ++(180+60:1.8) circle (0.15);
    \draw[fill] (6,0) ++(180+90:1.8) circle (0.15);
    \draw[fill] (6,0) ++(180+120:1.8) circle (0.15);
    \draw[fill] (6,0) ++(180+150:1.8) circle (0.15);
    
    \end{tikzpicture} 
\caption{A graph $G$ on $30$ vertices for which $\pi(G)+\partial_{\lfloor \frac{2D}{3}\rfloor} \approx 0.4$.}\label{fig:Wagner30Vertices}
\end{figure}
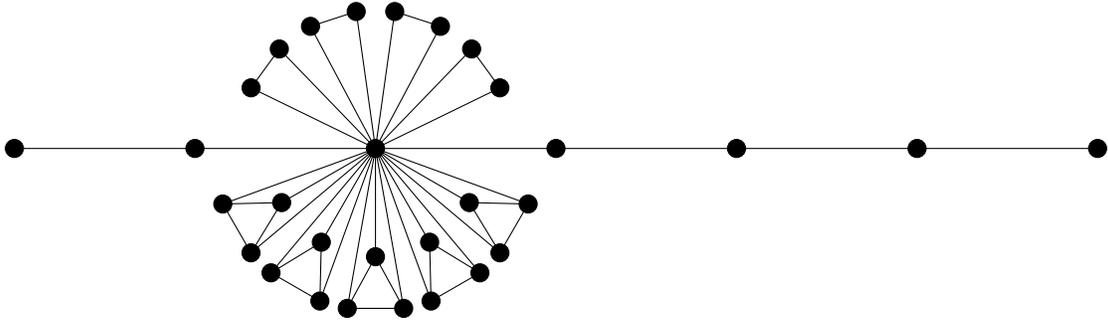

We now discuss a conjecture by Aouchiche and Hansen~\cite{AH16} that was disproved in this way. Let $G$ be a graph with vertices $\{v_1,v_2,\ldots,v_n\}$, let $C$ be its distance matrix (i.e., $C_{i,j}$ is equal to the distance $d(v_i,v_j)$ between $v_i$ and $v_j$ in $G$) and let $\partial_1 \geq \partial_2 \geq \ldots \geq \partial_n$ be the eigenvalues of this distance matrix $C$. Let the proximity $\pi(G) = \frac{1}{n-1} \cdot \min_{v \in V(G)} \sum_{w \in V(G)} d(v,w)$ and let $D = \max_{v,w \in V(G)}d(v,w)$ be the diameter of $G$. Aouchiche and Hansen~\cite{AH16} made the following conjecture:
\begin{conj}[\cite{AH16}]
\label{conj:AH}
Let $G$ be a connected graph on $n \geq 4$ vertices with diameter $D$, proximity $\pi(G)$ and eigenvalues of the distance matrix $\partial_1 \geq \partial_2 \geq \ldots \geq \partial_n$. Then:
$$\pi(G)+\partial_{\lfloor \frac{2D}{3}\rfloor}>0.$$
\end{conj}
This suggests to use $\pi(G)+\partial_{\lfloor \frac{2D}{3}\rfloor}$ as reward function (or $-1$ multiplied by this expression, depending on the precise implementation). After running the deep cross-entropy method for a few days using $n=30$ and the previously mentioned reward function, Wagner~\cite{W21} reports that the best graph that was found (see Fig.~\ref{fig:Wagner30Vertices}) has $\pi(G)+\partial_{\lfloor \frac{2D}{3}\rfloor} \approx 0.4$. While this graph is not directly a counterexample, its structure gives a good idea of what a potential counterexample could look like.

Wagner~\cite{W21} showed that this intuition is indeed correct: by taking a path on $13$ vertices and appending $190$ pendant vertices to the vertex at distance $4$ from one end on this path, one obtains a counterexample to Conjecture~\ref{conj:AH}.

In this subsection, we also feel compelled to draw the reader's attention to the recent successes of \textit{large language models}. These are neural networks that receive text-based input and also output text. This text can for example represent natural language or code in a programming language. Although we are not aware of any breakthrough results in graph theory being achieved by these models, they have led to new discoveries in disciplines similar to graph theory such as combinatorics (and in fact many other areas of mathematics) and we believe that it is just a matter of time before such models also lead to important new results in graph theory. We now give two examples of combinatorial results.

The \textit{cap set problem}~\cite{G19} is concerned with finding the largest set of vectors in $\mathbb{Z}^n_3$ such that no three vectors sum to zero. The authors of~\cite{REtAl24} proposed \texttt{FunSearch}: a procedure that combines evolutionary computation and large language models to write increasingly better code. Among other things, \texttt{FunSearch} generated code that outputs a set of $512$ vectors in $\mathbb{Z}^8_3$ such that no three vectors sum to zero, thereby improving the previously best known lower bound. Similarly, the authors of~\cite{NEtAl25} proposed \texttt{AlphaEvolve}, which substantially enhances \texttt{FunSearch} in several ways. For example, \texttt{AlphaEvolve} can generate code from any programming language whereas \texttt{FunSearch} is focused on \texttt{Python}. \texttt{AlphaEvolve} generated code that was able to output a step function that led to an improved upper bound for Erd{\H{o}}s's minimum overlap problem~\cite{E55}.

\section{New results by elementary means using computer-assisted graph theory}
\label{sec:newResults}
One aim of this survey is to promote the use of computer-assisted graph theory since we feel that it is currently underused in the literature. While we discussed several advanced techniques in the previous part of the survey, we now conclude by demonstrating that some gaps in the literature can be filled quite easily using a computer-assisted graph theory approach. All results in this section were obtained by very elementary means, requiring minimal effort, thanks to the common practice of researchers in computer-assisted graph theory to make their programs publicly available. This allows other researchers to benefit from advanced tools that often require years to develop and refine.\footnote{We emphasize the importance of properly citing such tools when they are used.}

\subsection{Small independence number of triangle-free graphs with given maximum degree}
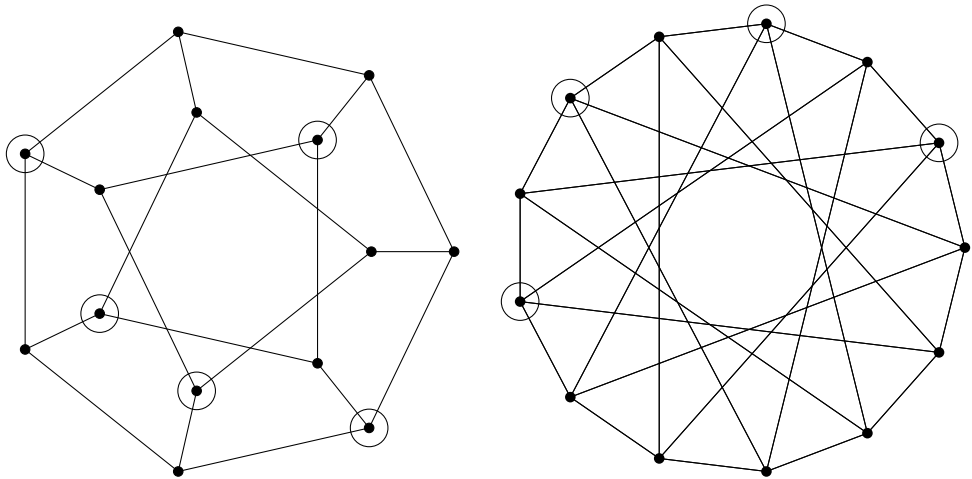
\begin{figure}[h!]
\begin{center}
\begin{tikzpicture}[scale=1.00]
        \foreach \x in {0,1,...,6}{
        \draw[fill] (\x*360/7:3.0) circle (1.8pt);
        \draw[fill] (\x*360/7:1.9) circle (1.8pt);
        \draw (\x*360/7:1.9)--(\x*360/7:3.0);
        \draw (\x*360/7+360/7:3.0)--(\x*360/7:3.0);
        \draw (\x*360/7+720/7:1.9)--(\x*360/7:1.9);
        }
        \draw[] (3*360/7+720/7:1.9) circle [radius=0.25];
        \draw[] (2*360/7+720/7:1.9) circle [radius=0.25];
        \draw[] (6*360/7+720/7:1.9) circle [radius=0.25];

        \draw[] (1*360/7+720/7:3.0) circle [radius=0.25];
        \draw[] (4*360/7+720/7:3.0) circle [radius=0.25];
        
        \end{tikzpicture} \quad
\begin{tikzpicture}
  \def\sides{13}
  \def\radius{3}

  \foreach \i in {1,...,\sides} {
    \fill ({360/\sides * \i}:\radius) circle (2pt);
  }

  \draw[] ({360/\sides * 1}:\radius) circle [radius=0.25];
  \draw[] ({360/\sides * 3}:\radius) circle [radius=0.25];
  \draw[] ({360/\sides * 5}:\radius) circle [radius=0.25];
  \draw[] ({360/\sides * 7}:\radius) circle [radius=0.25];
  
   \foreach \i in {1,...,\sides} {
    \draw ({360/\sides * (\i + 1)}:\radius) -- ({360/\sides * \i}:\radius);
    \draw ({360/\sides * (\i + 5)}:\radius) -- ({360/\sides * \i}:\radius);
    \draw ({360/\sides * (\i + 8)}:\radius) -- ({360/\sides * \i}:\radius);
    \draw ({360/\sides * (\i + 12)}:\radius) -- ({360/\sides * \i}:\radius);
  }
  
\end{tikzpicture}
\end{center}
\caption{Two graphs showing that $i(3) \leq \frac{5}{14} \approx 0.357$ (left) and $i(4) \leq \frac{4}{13} \approx 0.307$ (right). The circles indicate a maximum independent set.}\label{fig:i34} 
\end{figure}

\begin{figure}[h!]
\begin{center}
\begin{tikzpicture}
  \def\sides{20}
  \def\radius{3}

  \foreach \i in {1,...,\sides} {
    \fill ({360/\sides * \i}:\radius) circle (2pt);
  }

  \draw[] ({360/\sides * 5}:\radius) circle [radius=0.25];
  \draw[] ({360/\sides * 3}:\radius) circle [radius=0.25];
  \draw[] ({360/\sides * 7}:\radius) circle [radius=0.25];
  \draw[] ({360/\sides * -2}:\radius) circle [radius=0.25];
  \draw[] ({360/\sides * -4}:\radius) circle [radius=0.25];
  \draw[] ({360/\sides * -6}:\radius) circle [radius=0.25];
  
   \foreach \i in {1,...,\sides} {
    \draw ({360/\sides * (\i + 1)}:\radius) -- ({360/\sides * \i}:\radius);
    \draw ({360/\sides * (\i + 6)}:\radius) -- ({360/\sides * \i}:\radius);
    \draw ({360/\sides * (\i + 10)}:\radius) -- ({360/\sides * \i}:\radius);
  }
  
\end{tikzpicture} \quad
\begin{tikzpicture}
  \def\sides{28}
  \def\radius{3}

  \foreach \i in {1,...,\sides} {
    \fill ({360/\sides * \i}:\radius) circle (2pt);
  }
  
   \foreach \i in {1,8,15,22} {
    \draw ({360/\sides * (\i + 1)}:\radius) -- ({360/\sides * \i}:\radius);
    \draw ({360/\sides * (\i + 19)}:\radius) -- ({360/\sides * \i}:\radius);
    \draw ({360/\sides * (\i + 22)}:\radius) -- ({360/\sides * \i}:\radius);
    \draw ({360/\sides * (\i + 24)}:\radius) -- ({360/\sides * \i}:\radius);
    \draw ({360/\sides * (\i + 27)}:\radius) -- ({360/\sides * \i}:\radius);
  }

  \foreach \i in {2,16} {
    \draw ({360/\sides * (\i + 1)}:\radius) -- ({360/\sides * \i}:\radius);
    \draw ({360/\sides * (\i + 6)}:\radius) -- ({360/\sides * \i}:\radius);
    \draw ({360/\sides * (\i + 10)}:\radius) -- ({360/\sides * \i}:\radius);
    \draw ({360/\sides * (\i + 27)}:\radius) -- ({360/\sides * \i}:\radius);
  }

  \foreach \i in {3,17} {
    \draw ({360/\sides * (\i + 1)}:\radius) -- ({360/\sides * \i}:\radius);
    \draw ({360/\sides * (\i + 21)}:\radius) -- ({360/\sides * \i}:\radius);
    \draw ({360/\sides * (\i + 25)}:\radius) -- ({360/\sides * \i}:\radius);
    \draw ({360/\sides * (\i + 27)}:\radius) -- ({360/\sides * \i}:\radius);
  }

  \foreach \i in {4,6,18,20} {
    \draw ({360/\sides * (\i + 1)}:\radius) -- ({360/\sides * \i}:\radius);
    \draw ({360/\sides * (\i + 4)}:\radius) -- ({360/\sides * \i}:\radius);
    \draw ({360/\sides * (\i + 7)}:\radius) -- ({360/\sides * \i}:\radius);
    \draw ({360/\sides * (\i + 9)}:\radius) -- ({360/\sides * \i}:\radius);
    \draw ({360/\sides * (\i + 27)}:\radius) -- ({360/\sides * \i}:\radius);
  }

  \foreach \i in {5,19} {
    \draw ({360/\sides * (\i + 1)}:\radius) -- ({360/\sides * \i}:\radius);
    \draw ({360/\sides * (\i + 7)}:\radius) -- ({360/\sides * \i}:\radius);
    \draw ({360/\sides * (\i + 16)}:\radius) -- ({360/\sides * \i}:\radius);
    \draw ({360/\sides * (\i + 18)}:\radius) -- ({360/\sides * \i}:\radius);
    \draw ({360/\sides * (\i + 27)}:\radius) -- ({360/\sides * \i}:\radius);
  }

  \foreach \i in {7,21} {
    \draw ({360/\sides * (\i + 1)}:\radius) -- ({360/\sides * \i}:\radius);
    \draw ({360/\sides * (\i + 4)}:\radius) -- ({360/\sides * \i}:\radius);
    \draw ({360/\sides * (\i + 12)}:\radius) -- ({360/\sides * \i}:\radius);
    \draw ({360/\sides * (\i + 21)}:\radius) -- ({360/\sides * \i}:\radius);
    \draw ({360/\sides * (\i + 27)}:\radius) -- ({360/\sides * \i}:\radius);
  }

  \foreach \i in {9,23} {
    \draw ({360/\sides * (\i + 1)}:\radius) -- ({360/\sides * \i}:\radius);
    \draw ({360/\sides * (\i + 6)}:\radius) -- ({360/\sides * \i}:\radius);
    \draw ({360/\sides * (\i + 10)}:\radius) -- ({360/\sides * \i}:\radius);
    \draw ({360/\sides * (\i + 19)}:\radius) -- ({360/\sides * \i}:\radius);
    \draw ({360/\sides * (\i + 27)}:\radius) -- ({360/\sides * \i}:\radius);
  }

  \foreach \i in {10,24} {
    \draw ({360/\sides * (\i + 1)}:\radius) -- ({360/\sides * \i}:\radius);
    \draw ({360/\sides * (\i + 7)}:\radius) -- ({360/\sides * \i}:\radius);
    \draw ({360/\sides * (\i + 16)}:\radius) -- ({360/\sides * \i}:\radius);
    \draw ({360/\sides * (\i + 24)}:\radius) -- ({360/\sides * \i}:\radius);
    \draw ({360/\sides * (\i + 27)}:\radius) -- ({360/\sides * \i}:\radius);
  }

  \foreach \i in {11,25} {
    \draw ({360/\sides * (\i + 1)}:\radius) -- ({360/\sides * \i}:\radius);
    \draw ({360/\sides * (\i + 4)}:\radius) -- ({360/\sides * \i}:\radius);
    \draw ({360/\sides * (\i + 21)}:\radius) -- ({360/\sides * \i}:\radius);
    \draw ({360/\sides * (\i + 24)}:\radius) -- ({360/\sides * \i}:\radius);
    \draw ({360/\sides * (\i + 27)}:\radius) -- ({360/\sides * \i}:\radius);
  }

  \foreach \i in {12,26} {
    \draw ({360/\sides * (\i + 1)}:\radius) -- ({360/\sides * \i}:\radius);
    \draw ({360/\sides * (\i + 12)}:\radius) -- ({360/\sides * \i}:\radius);
    \draw ({360/\sides * (\i + 18)}:\radius) -- ({360/\sides * \i}:\radius);
    \draw ({360/\sides * (\i + 21)}:\radius) -- ({360/\sides * \i}:\radius);
    \draw ({360/\sides * (\i + 27)}:\radius) -- ({360/\sides * \i}:\radius);
  }
  
  \foreach \i in {13,27} {
    \draw ({360/\sides * (\i + 1)}:\radius) -- ({360/\sides * \i}:\radius);
    \draw ({360/\sides * (\i + 9)}:\radius) -- ({360/\sides * \i}:\radius);
    \draw ({360/\sides * (\i + 19)}:\radius) -- ({360/\sides * \i}:\radius);
    \draw ({360/\sides * (\i + 21)}:\radius) -- ({360/\sides * \i}:\radius);
    \draw ({360/\sides * (\i + 27)}:\radius) -- ({360/\sides * \i}:\radius);
  }

  \foreach \i in {14,28} {
    \draw ({360/\sides * (\i + 1)}:\radius) -- ({360/\sides * \i}:\radius);
    \draw ({360/\sides * (\i + 3)}:\radius) -- ({360/\sides * \i}:\radius);
    \draw ({360/\sides * (\i + 7)}:\radius) -- ({360/\sides * \i}:\radius);
    \draw ({360/\sides * (\i + 9)}:\radius) -- ({360/\sides * \i}:\radius);
    \draw ({360/\sides * (\i + 27)}:\radius) -- ({360/\sides * \i}:\radius);
  }

  \draw[] ({360/\sides * 1}:\radius) circle [radius=0.25];
  \draw[] ({360/\sides * 3}:\radius) circle [radius=0.25];
  \draw[] ({360/\sides * 5}:\radius) circle [radius=0.25];
  \draw[] ({360/\sides * 7}:\radius) circle [radius=0.25];
  \draw[] ({360/\sides * 9}:\radius) circle [radius=0.25];
  \draw[] ({360/\sides * 13}:\radius) circle [radius=0.25];
  \draw[] ({360/\sides * 16}:\radius) circle [radius=0.25];
  \draw[] ({360/\sides * 18}:\radius) circle [radius=0.25];

\end{tikzpicture}
\end{center}
\caption{Two graphs showing that $i(5) \leq \frac{6}{20}=0.3$ (left) and $i(5) \leq \frac{8}{28} \approx 0.286$ (right). The circles indicate a maximum independent set.}\label{fig:i5} 
\end{figure}
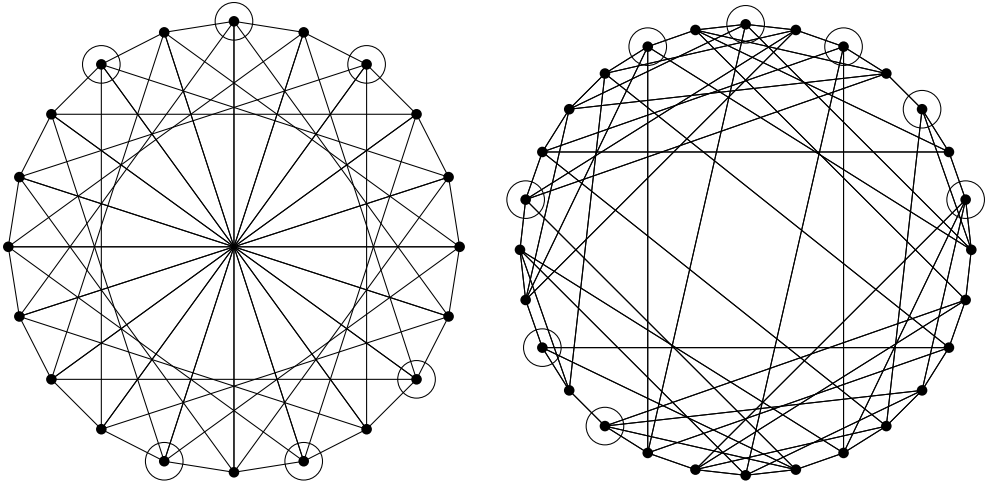

In this subsection, we are interested in the question how small the independence number $\alpha(G)$ (i.e., the size of a maximum independent set in $G$) can be for a triangle-free graph $G$ with maximum degree $\Delta$ compared to the order of $G$. More precisely, we are interested in determining the following: 
\[
\small
i(\Delta) := \inf\Bigl\{\frac{\alpha(G)}{n}~|~ G\text{ is a triangle-free graph on }n\text{ vertices with maximum degree at most }\Delta\Bigr\}.
\]
This quantity has received much attention in the literature. Staton~\cite{Staton79} proved that $i(3) \geq \frac{5}{14} \approx 0.357$, Jones~\cite{J84} proved that $i(4) \geq \frac{4}{13} \approx 0.307$ and Shearer~\cite{S91} proved that $i(5) \geq \frac{593}{2210} \approx 0.268$. For $i(3)$ and $i(4)$ it is known that these bounds are best possible due to the graphs shown in Fig.~\ref{fig:i34}, whereas for $i(5)$ the best upper bound is $i(5) \leq \frac{6}{20}=0.3$ due to Pirot and Sereni~\cite{PS21} (the corresponding graph is shown on the left of Fig.~\ref{fig:i5}).


However, it turns out that it is quite easy to achieve a better upper bound for $i(5)$ using the database House of Graphs~\cite{HOG} by searching for triangle-free graphs $G$ on $n$ vertices of maximum degree $5$ for which $\frac{\alpha(G)}{n}<0.3$. More precisely, there exists a minimal Ramsey graph~\cite{GR13} that was uploaded to the House of Graphs in 2012 (see the right side of Fig.~\ref{fig:i5} or \url{https://houseofgraphs.org/graphs/6540}) from which it follows that $i(5) \leq \frac{8}{28} \approx 0.286$. This perfectly illustrates the idea behind this database that many graphs that are important for one problem are also important for other problems.


\subsection{Planar graphical degree sequences}
A sequence of $n$ integers $\mathbf{d}=(d_1,d_2,\ldots,d_n)$ is called \textit{graphical} if there exists a graph having degree sequence $\mathbf{d}$. Likewise, a graphical sequence $\mathbf{d}$ is called \textit{planar graphical} if there exists a planar graph having degree sequence $\mathbf{d}$ (recall that we only deal with simple graphs). We will always assume that $d_1 \geq d_2 \geq \ldots \geq d_n$ and will use the shorthand notation $a^bc^d\ldots$ for a sequence consisting of $b$ repetitions of $a$, followed by $d$ repetitions of $c$ and so on. The Erd{\H{o}}s-Gallai theorem~\cite{EG60} states that $\mathbf{d}$ is graphical if and only if $\sum_{i=1}^n d_i$ is even and for every $k \in \{1,2,\ldots,n\}$ we have $\sum_{i=1}^k d_i \leq k(k-1)+\sum_{i=k+1}^n \min(d_i,k)$. The complete characterization of which degree sequences are planar graphical is harder than for graphical sequences and there are many partial results, but no complete characterization yet. 

Due to Euler's theorem, we have $\sum_{i=1}^n d_i \leq 6n-12$ for any planar graphical sequence. We call a planar graphical sequence \textit{maximal} if $\sum_{i=1}^n d_i = 6n-12$ and we call it a $k$-sequence if $d_1-d_n=k$. Schmeichel and Hakimi~\cite{SH77} established the following result:
\begin{thr}[\cite{SH77}]
Every graphical non-maximal $2$-sequence is planar graphical except for $4^52^1,5^{11}3^1, 5^53^3, 7^15^{15}, 6^{n-7}4^7$ for $n>7$, and possibly the following cases: $5^{13}3^1, 7^15^{17},7^35^{17}$.
\end{thr}

Hence, this leaves $3$ unresolved cases and Schmeichel and Hakimi made the following conjecture~\cite{SH77}:
\begin{conj}[\cite{SH77}]
The sequences $5^{13}3^1, 7^15^{17}$ and $7^35^{17}$ are not planar graphical.
\end{conj}

Fanelli~\cite{F81} was later able to show that the sequences $5^{13}3^1$ and $ 7^15^{17}$ are indeed not planar graphical, but was unable to solve the last missing case $7^35^{17}$. Up until now, this last case remained open (see for example~\cite{BBPRR24}). 

Using the graph generator \texttt{plantri}~\cite{BGGMTW05,BM05,BM07,VM24}, however, it is quite easy to show that $7^35^{17}$ is not planar graphical, thereby proving the last missing case of the Schmeichel-Hakimi conjecture. More precisely, generating all planar graphs on $20$ vertices would likely be too slow, but if there exists a planar graph $G$ with degree sequence $7^35^{17}$, then there exists a planar triangulation $G'$ on $20$ vertices such that $G$ can be obtained by removing an edge from $G'$ (since $3\cdot7+17\cdot5=106=6\cdot20-14$). Generating the planar triangulations using \texttt{plantri}, removing one edge in all possible ways and then checking the degree sequences yields the result in a matter of minutes.

\section*{Acknowledgements}

The author is supported by a Postdoctoral Fellowship of the Research Foundation Flanders (FWO) with grant number 1222524N. The author is grateful to Jan Goedgebeur and Brendan D. McKay for several discussions and suggestions that led to an improved presentation of this survey. 
\bibliographystyle{abbrv}
\bibliography{ref}

\end{document}